\documentclass{elsart}
\usepackage{amsmath}
\journal{Linear Algebra and its Applications}
\theoremstyle{remark}

\newcommand{\ra}{\rightarrow}
\newcommand{\la}{\leftarrow}
\newcommand{\nj}{$(0,1)$\ }
\newcommand{\dn}{$\mathcal{D}_n$\ }
\newcommand{\dnn}{$\mathcal{D}_n$}
\newcommand{\sn}{$\mathcal{S}_n$\ }
\newcommand{\snn}{$\mathcal{S}_n$}
\newcommand{\an}{$\mathcal{A}_n$\ }
\newcommand{\ann}{$\mathcal{A}_n$}
\newcommand{\ao}{$\mathcal{A}_8$\ }
\newcommand{\aon}{$\mathcal{A}_8$}
\DeclareMathOperator{\diag}{diag}
\DeclareMathOperator{\bord}{bord}
\DeclareMathOperator{\rang}{rank}
\DeclareMathOperator{\SNF}{SNF}
\newcommand{\hyph}{-\penalty0\hskip0pt\relax }
\newcommand{\rz}{\,}
\newcommand{\ry}{\hspace{0.7ex}}
\newcommand{\rg}{\,}
\newcommand{\rx}{\hspace{0.1ex}}
\newcommand{\rzz}{}
\newcommand{\br}{\makebox[1.16ex]{$\cdot$}}
\newcommand{\bn}{$\circ$}
\newcommand{\bj}{$\bullet$}
\newcommand{\bd}{$\star$}
\newcommand{\bc}{\bd}
\newcommand{\bs}{\bd}
\newcommand{\opt}{\mathrm{Optimize}}
\newcommand{\INPUT}{$\mathbf{Input:\ }$}
\newcommand{\OUTPUT}{$\mathbf{Output:\ }$}
\newcommand{\FOR}{$\mathbf{for\ }$}
\newcommand{\WHILE}{$\mathbf{while\ }$}
\newcommand{\DO}{$\mathbf{do\ }$}
\newcommand{\IF}{$\mathbf{if\ }$}
\newcommand{\THEN}{$\mathbf{then\ }$}
\newcommand{\ELSE}{$\mathbf{else\ }$}
\newcommand{\CONTINUE}{$\mathbf{continue;\ }$}

\begin{document}
\begin{frontmatter}
\title{Classification of small \nj matrices}
\author{Miodrag \v Zivkovi\'c}
\address{11000 Beograd\\Marka \v Celebonovi\'ca 61/15\\
Serbia and Montenegro}
\ead{ezivkovm@matf.bg.ac.yu}
\thanks{Supported by
 Ministry of Science and Technology of Serbia, Grant 1861}

\begin{abstract}
Denote by \an the set of square \nj matrices of order $n$.
The set \ann, $n\le8$, is partitioned into 
row/column permutation equivalence classes
enabling derivation of various facts by simple counting.
For example, the number of regular \nj
matrices of order $8$ is $10160459763342013440$. 
Let \dnn, \sn denote the set of absolute determinant values
and  Smith normal forms of matrices from \ann.
Denote by $a_n$ the smallest integer not in \dnn.
The sets
$\mathcal{D}_9$ and $\mathcal{S}_9$ are obtained;
especially, $a_9=103$.
The lower bounds for $a_n$, $10\le n\le 19$, 
(exceeding the known lower bound $a_n\ge 2f_{n-1}$,
where $f_n$ is $n$th Fibonacci number) are obtained.
Row/permutation equivalence classes of \an 
correspond to bipartite graphs with $n$ black and $n$
white vertices, and so the other applications 
of the classification are possible.
\end{abstract}
\end{frontmatter}

\section{Introduction}

Let \an denote the set of square \nj matrices of order $n$.
Hadamard maximum determinant problem is:
find the maximum determinant among the matrices in \ann.
In this paper
we consider a slightly more general problem:
determine the set 
$\mathcal{D}_n=\{|det A| \mid A\in \mathcal{A}_n$\}.

It is known~\cite{Willia46} that determinants of \nj matrices 
of order $n$ are related to
determinants of $\pm 1$ matrices of order $n+1$.
If $A$ is a $(0,1)$-matrix of order $n$, let $B=\Psi(A)$
be a $\pm 1$-matrix
 of order $n+1$ obtained from $A$
 by replacing its $0$ by $-1$, bordering with a row $-1$'s on the top,
 and a column of $1$'s on the right.
 Clearly, $\Psi$ is a one-to-one correspondence.
 By adding row $1$ of $B$ to each of the other rows of $B$,
 we see that
$\det B=2^n \det A$.


By the Hadamard inequality $|\det B|\le\sqrt{(n+1)^{n+1}},$ and
therefore for 
all 
$A\in\mathcal{A}_n$ 
$|\det A|\le2^{-n}\sqrt{(n+1)^{n+1}}$. 
The equality is attained if $B$ is
an Hadamard matrix, i.e. if $BB^T=(n+1)I_{n+1}$, where $T$ denotes
transposition, and $I_n$ is the unit matrix of order $n$; for
$n>2$ this implies $n=4k-1$. For upper bounds for determinants of
$A\in\mathcal{A}_n$ 
see for example~\cite{Neub97}.

Let $d_n$ denote the
largest element in \dnn, and let $a_n$ be the smallest integer not
in \dnn. Craigen~\cite{Craig} shows that the set \dn is
the interval $\{1,2,\ldots,d_n\}$ for $n\le 6$, but not
for $n=7$, because 
$a_8=41<d_8=56$; 
he suggests that $a_9=103$.

Some
interesting sequences, related to \nj matrices
are found in~\cite{Slo}: A003432 (the
sequence $d_n$), A013588 (the sequence $a_n$), A051752 ($c_n$, the
number of matrices in \an with the determinant $d_n$) and A055165
($m_n$, the number of regular matrices in \ann). A few first
members of these sequences are given in the following table.
The values of $a_9$, $c_8$, $c_9$ and $m_8$ seem to be new.

\begin{tabular}{||r||r|r|r|r||}
\hline\hline
      & A003432  &  A013588 & A051752      &  A055165  \\
$n$   &   $d_n$  & $a_n$    & $c_n$        &  $m_n$    \\
\hline\hline
$ 1$  &$   1$    & $   2$   & $1$          &$1$\\
$ 2$  &$   1$    & $   2$   & $3$          &$6$\\
$ 3$  &$   2$    & $   3$   & $3$          &$174$\\
$ 4$  &$   3$    & $   4$   & $60$         &$22560$           \\
$ 5$  &$   5$    & $   6$   & $3600$       &$12514320$        \\
$ 6$  &$   9$    & $  10$   & $529200$     &$28836612000$     \\
$ 7$  &$  32$    & $  19$   & $75600$      &$270345669985440$ \\
$ 8$  &$  56$    & $  41$   & $^*195955200$  &$^*10160459763342013440$\\
$ 9$  &$ 144$    & $^*103$  & $^*13716864000$& \\
$10$  &$ 320$    & $    $   &              &            \\
$11$  &$1458$    & $    $   &              &            \\
$12$  &$3645$    & $    $   &              &            \\
$13$  &$9477$    & $    $   &              &            \\
\hline\hline
\end{tabular}


In this paper, which is a continuation of~\cite{zivkov01}, 
the matrices in \ann, $n\le8$, are
partitioned into row/column permutation equivalence classes,
enabling the
classification by ADV, and more precisely --- by
 SNF (see section~\ref{razvrst}).
Let
$\mathcal{S}_n$ denote the set of
 SNF's of matrices in
$\mathcal{A}_n$.
In section~\ref{snf9} the sets
$\mathcal{D}_9$ and
$\mathcal{S}_9$ are determined.
In section~\ref{donjegr} the lower bounds
for $a_n$, $10\le n\le 19$ are obtained;
 $c_n$, $n\le 9$,  are obtained in section~\ref{brojevi}.

We introduce now some notation.
If
$A=[a_{ij}]$ and
$B=[b_{ij}]$
are matrices of the same dimension $m\times n$,
we say that $A<B$ if
$A$ is \emph{lexicographically less than} $B$, i.e. if for some pair of indices
$(i,j)$
the first $i-1$ rows of $A$ and $B$ are equal,
the first $j-1$ elements in the $i$th row of $A$ and $B$ are equal,
 and $a_{ij}<b_{ij}$.
For example,
\[
\left[\begin{array}{rr}
  1 & 0 \\
  1 & 0
\end{array}\right] <
\left[\begin{array}{rr}
  1 & 0 \\
  1 & 1
\end{array}\right].
\]

The smallest matrix in a set 
$\mathcal{A}$
is the representative of 
$\mathcal{A}$.

Denote by
$P_{i,j}$
the permutation matrix obtained from $I_n$
by exchanging the
$i$th and $j$th row.

The matrices $A,B\in\mathcal{A}_n$ are \emph{equivalent}~\cite{gantmaher},
$A\sim B$, if $B$ is be obtained from $A$ by a sequence of
elementary row/column operations of the following types: exchange
of two rows/columns, multiplication of a row/column by $-1$, and
addition/subtraction of a row/column to/from another row/column. Let
$\SNF(A)$ denote the SNF of $A$. It is known that $A\sim B$ is
equivalent to $\SNF(A)=\SNF(B)$ 
(in~\cite{gantmaher} this statement is proved for 
polynomial matrices).

The SNF $\diag(d_1,d_2,\ldots,d_n)$ is written simply as a vector
$(d_1,d_2,\ldots,d_n)$. If diagonal elements of SNF are repeated, we use
the shortened exponential notation. For example,
$(1^3,2,0)$ is short
$(1,1,1,2,0)$.
If
$s\in \mathcal{S}_n$, then
we also say that the SNF-class $s$ is the set
$\{A\in\mathcal{A} \mid \SNF(A)=s\}$.

Let $J_n$ denote the square matrix of order $n$
with all elements equal to one.

\section{Classification of \nj matrices of order $8$ or less}\label{razvrst}

The set
\dn could be obtained by computing
determinants of all $A\in \mathcal{A}_n$.
A better approach is to group matrices with the same determinant,
and then to compute
the determinant of only one matrix in each group.
It is useful to classify
\an
into subsets with constant absolute determinant value(ADV),
or into even smaller subsets with constant SNF.
We now review some such partitions of $\mathcal{A}_n$.

Let
$\Pi_r$ denote the group of row permutations of matrices
from \ann.
Permutations from $\Pi_r$ preserve ADV.

The representative of the matrix $A$ orbit
is obtained from
$A$ by sorting its rows into a nondecreasing sequence.
Rows of $A$ correspond to binary numbers  
less than $N=2^n$.
Therefore, the number of orbits of
$\Pi_r$ in \an is equal to 
$\binom{N+n-1}{n-1}$, i.e.
the number of nondecreasing sequences of length $n$ 
from $\{0,1,\ldots,N-1\}$.
Let $\Pi$ denote the group of row and column permutations; 
$\Pi$ also preserves ADV. The group $\Pi$ induces
an equivalence relation $\pi$ over \ann.
We say that matrices $A$ and $B$
are permutationally  equivalent, $A\sim_\pi B$, if they are in the
same orbit of $\Pi$. Let $A_\pi$ denote
the representative of the matrix
$A$ equivalence class
($\pi$\hyph class; we say shorter that $A_\pi$  is a
$\pi$\hyph representative of $A$).

\begin{exmp}\label{ex1}
The $\pi$\hyph representative of
\[
\left[\begin{array}{rrr}
  1 & 0 & 1 \\
  1 & 1 & 0 \\
  1 & 0 & 0
\end{array}\right]
\]
is the matrix
\[
\left[\begin{array}{rrr}
  0 & 0 & 1 \\
  0 & 1 & 1 \\
  1 & 0 & 1
\end{array}\right],
\]
the smallest of all $36$ permutationally equivalent matrices.
\end{exmp}

Let $\mathcal{A}_n^\pi$ denote the set of $\pi$-representatives in \ann.
In~\cite{Harary} it is shown that the number of
$\pi$\hyph classes in $\mathcal{A}_n$ is given by:
\begin{equation}\label{pclass}
\vert \mathcal{A}_n^\pi\vert=
\sum_{i_1+2i_2+\ldots+ni_n=n} 
\ \sum_{j_1+2j_2+\ldots+nj_n=n} 
C(i) C(j)
\exp_2\sum_{r,s=1}^n i_r j_s 2^{(r,s)},
\end{equation}
where the summation is over all vectors
$i=(i_1,i_2,\ldots,i_n)$,
$j=(j_1,j_2,\ldots,j_n)$,
and
\[C(i)=n!/(1^{i_1}i_1!\ldots n^{i_n}i_n!)\]
is the number $n$\hyph permutations with $i_r$ cycles of length $r$,
$r=1,2,\ldots,n$; 
$(r,s)$ denotes GCD of integers $r$, $s$. 
The values $\vert \mathcal{A}_n^\pi\vert$ are
listed in Table~\ref{ubas}; they are easily  computed for quite a
large $n$ using, for example, UBASIC~\cite{UB}. 
It is seen that $p_n$ is close to
$2^{n^2}/(n!)^2$ for $n\le15$. An effective algorithm to generate
the representative $A_\pi$ of a given matrix $A$ 
(section~\ref{spipredst}) simplifies the
classification of matrices, because it enables to deal with 
the small subset $\mathcal{A}_n^\pi$ of \ann.

\begin{table}
\caption{The number of permutationally
nonequivalent matrices in \ann, $n\le 15$.}
\label{ubas}
{\small
\begin{tabular}{||r|r|r||}
\hline\hline
$n$ & $(2^{n^2}/n!^2)/|\mathcal{A}_n^\pi |$ &$|\mathcal{A}_n^\pi |$    \rule[-2mm]{0cm}{6mm}   \\
\hline\hline
  1&  1.00000&                                             2\\\hline
  2&  0.57143&                                             7\\\hline
  3&  0.39506&                                            36\\\hline
  4&  0.35892&                                           317\\\hline
  5&  0.41433&                                          5624\\\hline
  6&  0.52685&                                        251610\\\hline
  7&  0.65875&                                      33642660\\\hline
  8&  0.77266&                                   14685630688\\\hline
  9&  0.85533&                                21467043671008\\\hline
 10&  0.91045&                            105735224248507784\\\hline
 11&  0.94565&                        1764356230257807614296\\\hline
 12&  0.96754&                   100455994644460412263071692\\\hline
 13&  0.98088&              19674097197480928600253198363072\\\hline
 14&  0.98886&        13363679231028322645152300040033513414\\\hline
 15&  0.99358&  31735555932041230032311939400670284689732948\\\hline
\hline\hline
\end{tabular}

}
\end{table}

\subsection{Matrix extension}

In order to classify
matrices in \an by ADV values, one has to select carefully
the order by which determinants are computed.
It is natural to start from
matrices of order $n-1$,
and then to extend them by one row and one column of
ones and zeros in each possible way.
For an arbitrary $B\in\mathcal{A}_{n-1}$, let
$\bord(B)$
denote the subset of  \ann, 
containing matrices with the upper left
minor equal to $B$. We say that the matrices in
$\bord(B)$ are \emph{obtained by extending}
$B$;
if $A\in \bord(B)$, then $A$ is an \emph{extension} of $B$.

The calculation of determinants of all matrices in $\bord(B)$
is an easy task.
If $A\in\bord(B)$, then $A$ is of the form
\begin{equation}\label{bordA}
A=
\left[\begin{array}{rr}
  B & y \\
  x & b
\end{array}\right],
\end{equation}
where $x=[x_1\ x_2\ \ldots\ x_{n-1}]$ and $y=[y_1\ y_2\ \ldots\
y_{n-1}]^T$. Then~\cite{Willia46}
\begin{equation}\label{detobr}
\det A=b\det B -
\sum_{i=1}^{n-1}
\sum_{j=1}^{n-1}
x_i y_j \det B_{ij},
\end{equation}
where $B_{ij}$
is the cofactor of $B$, corresponding to $a_{ij}$.

Obviously,
\[
\mathcal{A}_{n}=\{ A \mid
(B,x,y,b)\in \mathcal{A}_{n-1}\times\{0,1\}^{n-1}\times\{0,1\}^{n-1}\times\{0,1\}\}.
\]
If we precompute cofactors
$B_{ij}$,
then determinant of each matrix from $\bord(B)$
is computed by only one addition:
for the fixed $x$, the column $y$ might
traverse the set of possible values via
a Gray code
(so that in the sequence of $y$'s each two subsequent vectors
differ in exactly one position).

Williamson~\cite{Willia46} noted that it is  
enough to let $B$ cross the set of
$\pi$\hyph representatives in 
$\mathcal{A}_{n-1}$.      
Let $\bord_\pi(B)$ denote the set of
$\pi$\hyph representatives of matrices in $\bord(B)$.

\begin{lem}\label{pobrub}
If
$B\sim_\pi B'$
then
$\bord_\pi(B)=\bord_\pi(B')$.
\end{lem}

\begin{pf}
Let $A\in\bord_\pi(B)$.
If the row/column permutations, transforming
$B$ into $B'$, are applied to the first
$n-1$ rows/columns of $A$,
then the
matrix with the upper left minor equal to $B'$ is obtained.
Therefore, the
matrix permutationally equivalent to $A$ is obtained
by extending $B'$,
meaning that
$A$ is permutationally equivalent to a matrix from $\bord(B')$, i.e.
$A\in \bord_\pi(B')$. Analogously,
$\bord_\pi(B')\subseteq \bord_\pi(B)$, and so
$\bord_\pi(B') = \bord_\pi(B)$.
\qed \end{pf}

Not only determinants, but also SNF's of
matrices in $\bord(B)$
can be efficiently computed.
The preprocessing step is to compute
$D=\SNF(B)=\diag(d_1,d_2,\ldots,d_n)$,
and the matrices
$P$, $Q$,
such that
$PBQ=D$,
$|\det P|=|\det Q|=1$.
In order to determine
$\SNF(A)$ for an arbitrary
 $A\in\bord(B)$ of the form~\eqref{bordA},
we use the identity
\begin{equation}\label{snfobr}
\left[\begin{array}{rr}
  P & 0 \\
  0 & 1
\end{array}\right]
\left[\begin{array}{rr}
  B & y \\
  x & b
\end{array}\right]
\left[\begin{array}{rr}
  Q & 0 \\
  0 & 1
\end{array}\right]=
\left[\begin{array}{rr}
  D & Py \\
  xQ & b
\end{array}\right].
\end{equation}
Denote
$xQ=[a_1\ a_2\ \ldots\ a_n]$,
$Py=[c_1\ c_2\ \ldots\ c_n]^T$.
Suppose $d_1=d_2=\cdots=d_k=1$, for some $k$, $1\le k\le n$.
Transforming the matrix from the righthand side of~\eqref{snfobr}
by subtracting the row
$i$ multiplied by $c_i$ from the row $n$, $1\le i\le k$,
 and then subtracting the column $i$
 multiplied by $c_i$ from the column $n$, $1\le i\le k$,
 we derive that $A$ is equivalent to
\begin{equation}\label{zz}
\left[\begin{array}{rrr|rrr|c}
  1   &\ldots &   0 &     0& \ldots&   0  &   0   \\
\vdots&\ddots&\vdots&\vdots&       &\vdots&\vdots \\
  0   &\ldots &   1 &     0& \ldots&   0  &   0   \\
\hline
  0   &\ldots &   0 &d_{k+1}&\ldots&   0  &c_{k+1}\\
\vdots&      &\vdots&\vdots&\ddots &\vdots&\vdots \\
  0   &\ldots &   0 &     0& \ldots&d_n   & c_n   \\
\hline
  0   &\ldots &   0 &a_{k+1}&\ldots&a_n   &b-\sum_{i=1}^k a_i c_i \\
\end{array}\right]
.
\end{equation}
Hence, $\SNF(A)$ determination is reduced
to determination
of SNF of a matrix of order $n-k$.
The special cases when
 $k\ge n-1$ are extremely simple, and they are not rare at all,
 because the corresponding SNF-classes are among the largest ones
 (at least for $n\le 9$).
More generally,
one can reduce $a_i$, $c_i$ modulo $d_i$,
$1\le i\le \rang B$.

\subsection{$\Phi$\hyph extension}

Following Williamson~\cite{Willia46}, the approach based on
extending  $\pi$\hyph representatives only, can be further improved.

For an arbitrary
$A\in\mathcal{A}_n$
let
$A'=X_i A$
denote the
matrix with the $i$th row equal to the
$i$th row of $A$, and with the row
$j\neq i$ equal to the coordinatewise modulo two sum
of $j$th and $i$th row of $A$.
Equivalently, $A'=RAS$, where
$R$ is the matrix obtained from $I_n$ by subtracting $i$th
row from the others, and
then by multiplying
 $i$th row by $-1$;
$S$ is the matrix obtained from $I_n$ by 
changing sign of columns corresponding to ones 
 in the $i$th row of $A$. A third equivalent
definition of $X_i$~\cite{Willia46} can be stated as follows: in
the $\ \pm1$ matrix $B=\Psi(A)$ of order $n+1$, the rows $1$ and
$(i+1)$ are exchanged, then the first row is ''normalized''
to all ones by changing signs of appropriate
columns. By applying
$\Psi^{-1}$, the matrix $A'$ is obtained. Therefore, application
of $X_i$ to $A$ corresponds to a special row permutation in
$\Psi(A)$ (followed by scaling).
It is natural to denote the identity transformation by
$X_0$, $X_0 A=A$.

The transformation $X_i$ also preserves ADV.
The composition of arbitrary two transformations
$X_i$, $X_j$
is equivalent to only one:
\[
X_i (X_j A) =
\left\{
 \begin{array}{ll}
P_{i,j} (X_i A),    & \text{if } 1\le i,j\le n \text{ if } i\neq j\\
 A,               &  \text{if } 1\le i=j \le n
 \end{array}
\right. .
\]
Let
$\Phi_r$ denote the set of $(n+1)!$ transforms of the form
$PX_i$, $0\le i\le n$,
where $P$ is an arbitrary permutation matrix.

\begin{thm}
The set $\Phi_r$ is a transformation group of \ann.
\end{thm}

\begin{pf}
We have
\[
X_i P A = P X_{p_i} A,
\]
where $p_i$ is the index of the row of $A$,
which is moved to the position $i$ after the left multiplication by
$P$.
Let $P_1$ and $P_2$ be the two permutation  matrices and
let $p_j$ be the position to which $P_1$ moves the row
$j$ after the left multiplication. Then
\[
P_2X_jP_1X_i =
P_2P_1X_{p_j}X_i =
\left\{
 \begin{array}{ll}
   P_2 P_1, & p_j=i \\
   P_2 P_1 P_{p_j,i} X_{p_j}, & p_j\neq i \\
 \end{array}
\right. .
\]
If $P_1=P$ is an arbitrary permutation matrix,
$1\le i\le n$,
$P_2=P^{-1}$,
 and
$p_j=i$, then
\[
 (P X_i)^{-1} = P^{-1} X_{j}.\qed
\]
 \end{pf}

Clearly, each orbit of
$\Phi_r$ contains at most $n+1$ orbits of $\Pi_r$.

The corresponding transformation $A X_j$
over the columns of $A$
(coordinatewise addition modulo two of the column $i$
to all other columns)
is defined by
$A X_j= (X_j A^T)^T$.
Let
$\Phi_c$
denote the group
generated by column permutations and
column transformations $(\cdot)X_i$.

Let $\Phi$ be the group generated by the elements of
groups
$\Phi_r$ and
$\Phi_c$; it also preserves ADV and 
its size is $(n+1)!^2$.
Matrices $A$ and $A'$ are said to be
$\phi$\hyph equivalent, $A\sim_\phi A'$, if they belong to the same orbit of
 $\Phi$.
 Equivalently, $A\sim_\phi A'$ if and only if there exist row and column
permutations
$P$, $Q$, and row and column transformations $X_i$, $X_j$,
such that
$A=P X_i A' X_j Q$.
For an arbitrary
$A\in\mathcal{A}_n$ let $A_\phi$ denote the
$\phi$\hyph representative of $A$;
$\phi$\hyph class of $A$ is the orbit of $\Phi$ containing $A$.

Let $\bord_\phi(B)$ denote the set of $\phi$\hyph representatives of
matrices in $\bord(B)$. Williams~\cite{Willia46} noted that
 $\Phi$ and $\Pi$ have similar properties:
 in order to obtain the set $\mathcal{A}_{n}^\phi$ of all
$\phi$\hyph representatives in \ann,
it is enough to extend
$\phi$\hyph representatives in $\mathcal{A}_{n-1}$.

\begin{lem}
If
$B\sim_\phi B'$,
then
$\bord_\phi(B)=\bord_\phi(B')$.
\end{lem}

\begin{pf}
If $B$ and $B'$ are $\phi$\hyph equivalent,
then there exist
$g\in\Phi$, transforming $B$ into $B'$.
Suppose
$A\in \bord_\phi(B)$.
Then
there exists a matrix $A'\in\bord(B)$,
$A'\sim_\phi A$.
By applying $g$ to upper left minor of $A$,
the matrix
$A''\sim_\phi A'$,
$A''\in\bord(B')$
is obtained.
Therefore,
$A\sim_\phi A''$,
and
$A\in \bord(B')$. Because $A$ is a $\phi$\hyph representative,
we obtain $A\in \bord_\phi (B')$, implying
$\bord_\phi(B)\subseteq \bord_\phi(B')$.
Analogously,
$\bord_\phi(B')\subseteq \bord_\phi(B)$, and hence
$\bord_\phi(B) = \bord_\phi(B')$.
\qed \end{pf}

\subsection{Effective determination of $\pi$\hyph representatives}
\label{spipredst}

The classification of  matrices in \an by extending matrices from
$\mathcal{A}_{n-1}^\phi$
must be accompanied
by an effective procedure to determine
$A_\pi$ and
$A_\phi$
for an arbitrary $A\in\mathcal{A}_n$.

The matrix $A_\pi$ is the smallest
among the family of at most $n!$ matrices
obtained by sorting rows of all the column permutations of $A$.
Search is performed more efficiently
by a branch-and-bound algorithm.
If we know the first $i$ rows
of $A_\pi$
(i.e. the row and column permutations $P$, $Q$ such that
the first $i$ rows of
$PAQ$ are minimal),
then the next row of
$A_\pi$
is a smallest column permutation
(only permutations
preserving the first $i$ rows of $PAQ$ are considered)
of some of the remaining rows of $PAQ$.

\begin{alg}\label{pipredst}
Branch-and-bound algorithm to determine
$A_\pi$ given $A\in\mathcal{A}_n$.
\begin{tabbing}
 \quad \= \quad \= \quad \= \quad  \kill
  \INPUT  $A\in\mathcal{A}_n$  \\
  \OUTPUT  $A_\pi$; the permutation matrices $P$, $Q$,
 such that $PAQ=A_\pi$;  \\
  \>   $count$ -- the number of pairs $(P,Q)$, such that
$PAQ=A_\pi$; \\
    $P^{(0)}\la I_n$; $Q^{(0)}\la I_n$; $A_\pi\la J_n$;  \\
    $i\la 0$;  \\
    $count\la 0$;  \\
    $\opt(i)$;  \\
      \\
\end{tabbing}

\begin{tabbing}
 \quad \= \quad \= \quad \= \quad  \kill
    \{Continuation of the search for $A_\pi$
starting from the row $i$ of $P^{(i-1)}AQ^{(i-1)}$,\} \\
    \{i.e. when the first $i-1$ rows are already chosen and permuted\} \\
    $\opt(i)$  \\
    Generate the minimal set of boundaries
$\Sigma^{(i)}= (s^{(i)}_0=0, s^{(i)}_1,\ldots, s^{(i)}_{k_i}=n)$  \\
  \>   between adjacent columns of $P^{(i-1)}AQ^{(i-1)}$, such that the $(i-1)$-prefixes \\
  \>  of columns from $s^{(i)}_{j-1}+1$ to $s^{(i)}_j$ are mutually equal, $1\le j\le k_i$;  \\
  \FOR  $j=i$ to $n$  \DO \\
  \>   $v_{jl}\la\sum_{r=s_{l-1}+1}^{s_l} \left( P^{(i-1)}AQ^{(i-1)} \right)_{j,r}$,
          $1\le l\le k_i$ ; \{the number of $1$'s \}  \\
  \> \>   \{in positions from $s^{(i)}_{j-1}+1$ to $s^{(i)}_j$ in the $j$th row of $P^{(i-1)}AQ^{(i-1)}$\}
\end{tabbing}

\begin{tabbing}
 \quad \= \quad \= \quad \= \quad  \kill
    $L^{(i)} \la$ the list of indices of largest vectors $v_j=(v_{j1},v_{j2},\ldots,v_{j{k_i}})$,
    $i\le j\le n$;  \\
  \FOR  all $j\in L^{(i)}$  \DO   \{     the candidates for the $i$th row of $A_\pi$\} \\
  \>   $P^{(i)}\la P_{i,j} P^{(i-1)}$;  \{    exchange the rows $i$ and $j$\}  \\
  \>   compute $Q^{(i)}$ from $Q^{(i-1)}$, so that all $1$'s in the part of the row $i$
  from  \\
  \> \>   $s^{(i)}_{l-1}+1$ to $s^{(l)}_i$  are moved to the right end of the part, $1\le l\le k_i$; \\
  \> \>   \{hence preserving the first $i-1$ rows of $P^{(i)}AQ^{(i)}$\} \\
  \>    compare the $i$th row of $P^{(i)}AQ^{(i)}$ to the $i$th row of $A_\pi$: \\
  \>  \IF  the $i$th row of $P^{(i)}AQ^{(i)}$ is less  \THEN\\
  \> \>    copy the first $i$ rows from $P^{(i)}AQ^{(i)}$ into $A_\pi$; \\
  \> \>    fill with ones the rest of $A_\pi$;\\
  \> \>  \IF  $i=n$  \THEN   $P\la P^{(i)}$; $Q\la Q^{(i)}$; $count\la 1$;
 \ELSE    $\opt(i+1)$; \\
  \>  \ELSE\IF  the $i$th row of $P^{(i)}AQ^{(i)}$ greater  \THEN\\
  \> \>  \CONTINUE   \{bound step: try the next row index from $L^{(i)}$\} \\
  \>  \ELSE\\
   \> \>  \IF  $i=n$  \THEN    $count\la count+1$;  \ELSE   $\opt(i+1)$; \\
\end{tabbing}
\end{alg}

\begin{exmp}\label{ex2}
Algorithm~\ref{pipredst}, applied to
the matrix from
Example~\ref{ex1},
gives the same $\pi$\hyph representative as obtained by trivial algorithm,
see Figure~\ref{ppipredst}.
\end{exmp}

\setlength{\unitlength}{0.9cm}
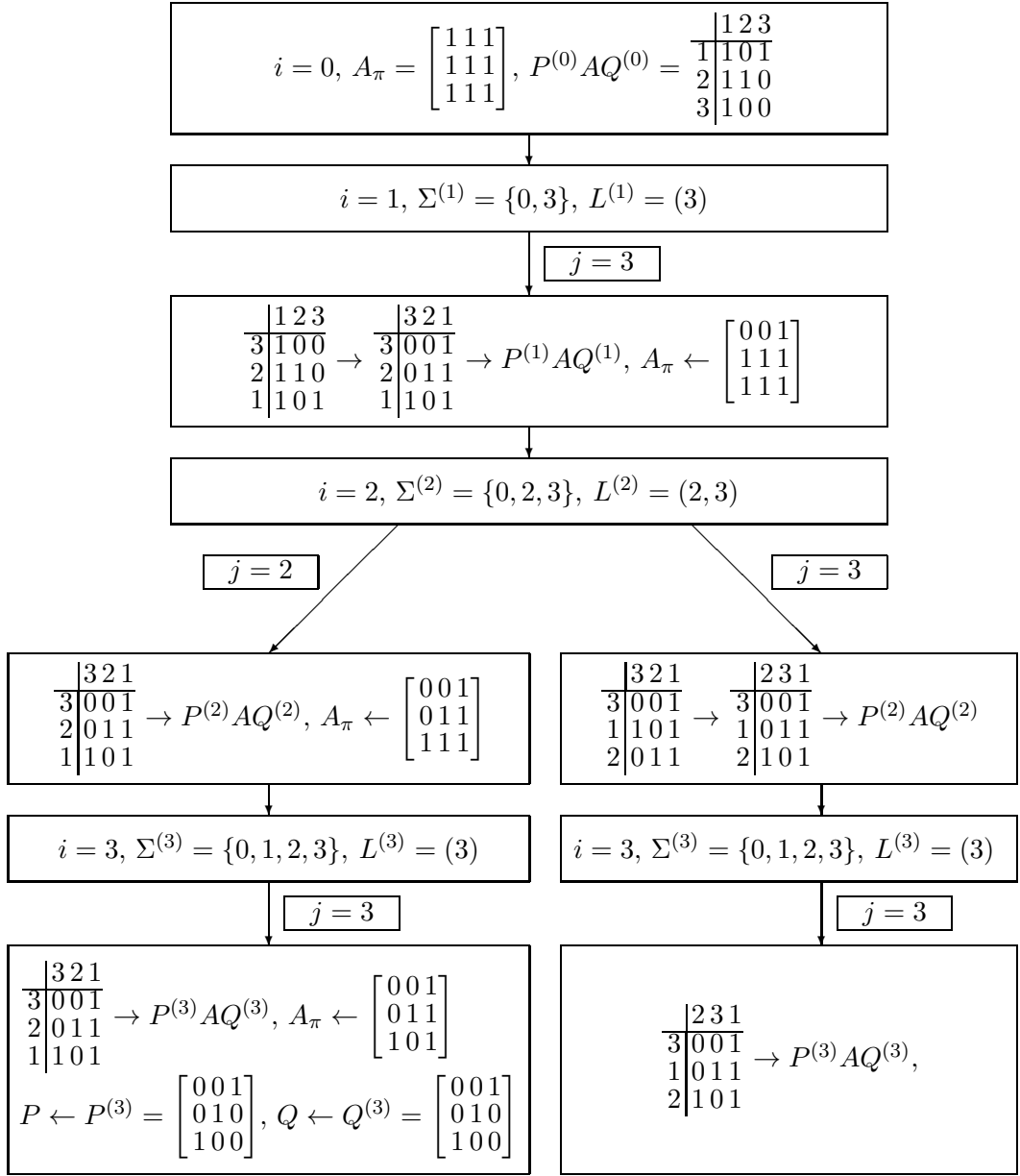
\begin{figure}
\renewcommand{\arraystretch}{0.8}
\begin{picture}(14,18)
\footnotesize
\put(2,16){\framebox(11,2){
$i=0$,
$A_\pi=
\left[\begin{array}{@{\rz}r@{\rz}r@{\rz}r@{\rzz}}
  1 & 1 & 1 \\
  1 & 1 & 1 \\
  1 & 1 & 1
\end{array}\right]$,
$P^{(0)}AQ^{(0)}=
\left.\begin{array}{r|r@{\rz}r@{\rz}r@{\rzz}}
   & 1 & 2 & 3 \\
   \hline
 1 & 1 & 0 & 1 \\
 2 & 1 & 1 & 0 \\
 3 & 1 & 0 & 0 \\
\end{array}\right.$
}}
\put(2,14.5){\framebox(11,1){
$i=1$, $\Sigma^{(1)}=\{0,3\}$, $L^{(1)}=(3)$
}}
\put(2,11.5){\framebox(11,2){
$\left.\begin{array}{r|r@{\rz}r@{\rz}r@{\rzz}}
   & 1 & 2 & 3 \\
   \hline
 3 & 1 & 0 & 0 \\
 2 & 1 & 1 & 0 \\
 1 & 1 & 0 & 1 \\
\end{array}\right. \ra%
\left.\begin{array}{r|r@{\rz}r@{\rz}r@{\rzz}}
   & 3 & 2 & 1 \\
   \hline
 3 & 0 & 0 & 1 \\
 2 & 0 & 1 & 1 \\
 1 & 1 & 0 & 1 \\
\end{array}\right. \ra P^{(1)}AQ^{(1)}$,
$A_\pi\la
\left[\begin{array}{@{\rz}r@{\rz}r@{\rz}r@{\rzz}}
  0 & 0 & 1 \\
  1 & 1 & 1 \\
  1 & 1 & 1
\end{array}\right]$
}}
\put(2,10){\framebox(11,1){$i=2$, $\Sigma^{(2)}=\{0,2,3\}$, $L^{(2)}=(2,3)$}}
\put(-0.5,6){\framebox(8,2){
$\left.\begin{array}{r|r@{\rz}r@{\rz}r@{\rzz}}
   & 3 & 2 & 1 \\
   \hline
 3 & 0 & 0 & 1 \\
 2 & 0 & 1 & 1 \\
 1 & 1 & 0 & 1 \\
\end{array}\right. \ra P^{(2)}AQ^{(2)}$,
$A_\pi\la
\left[\begin{array}{@{\rz}r@{\rz}r@{\rz}r@{\rzz}}
  0 & 0 & 1 \\
  0 & 1 & 1 \\
  1 & 1 & 1
\end{array}\right]$
}}
\put(-0.5,4.5){\framebox(8,1){
$i=3$, $\Sigma^{(3)}=\{0,1,2,3\}$, $L^{(3)}=(3)$
}}
\put(-0.5,0){\framebox(8,3.5){\,\,\,\,\,
\parbox{8cm}{\quad
$\left.\begin{array}{r|r@{\rz}r@{\rz}r@{\rzz}}
   & 3 & 2 & 1 \\
   \hline
 3 & 0 & 0 & 1 \\
 2 & 0 & 1 & 1 \\
 1 & 1 & 0 & 1 \\
\end{array}\right. \ra P^{(3)}AQ^{(3)}$,
$A_\pi\la
\left[\begin{array}{@{\rz}r@{\rz}r@{\rz}r@{\rzz}}
  0 & 0 & 1 \\
  0 & 1 & 1 \\
  1 & 0 & 1
\end{array}\right]$\\
\rule[-2mm]{0cm}{6mm}\quad{
$P\la P^{(3)}=
\left[\begin{array}{@{\rz}r@{\rz}r@{\rz}r@{\rzz}}
  0 & 0 & 1 \\
  0 & 1 & 0 \\
  1 & 0 & 0
\end{array}\right]$,
$Q\la Q^{(3)}=
\left[\begin{array}{@{\rz}r@{\rz}r@{\rz}r@{\rzz}}
  0 & 0 & 1 \\
  0 & 1 & 0 \\
  1 & 0 & 0
\end{array}\right]$
}
}
}}
\put(8,6){\framebox(7,2){
$\left.\begin{array}{r|r@{\rz}r@{\rz}r@{\rzz}}
   & 3 & 2 & 1 \\
   \hline
 3 & 0 & 0 & 1 \\
 1 & 1 & 0 & 1 \\
 2 & 0 & 1 & 1 \\
\end{array}\right. \ra
\left.\begin{array}{r|r@{\rz}r@{\rz}r@{\rzz}}
   & 2 & 3 & 1 \\
   \hline
 3 & 0 & 0 & 1 \\
 1 & 0 & 1 & 1 \\
 2 & 1 & 0 & 1 \\
\end{array}\right.\ra
P^{(2)}AQ^{(2)}$
}}
\put(8,4.5){\framebox(7,1){
$i=3$, $\Sigma^{(3)}=\{0,1,2,3\}$, $L^{(3)}=(3)$
}}
\put(8,0){\framebox(7,3.5){
$\left.\begin{array}{r|r@{\rz}r@{\rz}r@{\rzz}}
   & 2 & 3 & 1 \\
   \hline
 3 & 0 & 0 & 1 \\
 1 & 0 & 1 & 1 \\
 2 & 1 & 0 & 1 \\
\end{array}\right.\ra
P^{(3)}AQ^{(3)}$,
}}
\put(7.50,16.0){\vector(0,-1){0.5}}
\put(7.50,14.5){\vector(0,-1){1.}}
\put(7.75,13.75){\framebox(1.75,0.5){$j=3$}}
\put(7.5,11.5){\vector(0,-1){0.5}}
\put(2.500,9.){\framebox(1.75,0.5){$j=2$}}
\put(11.25,9.){\framebox(1.75,0.5){$j=3$}}
\put(5.5,10.){\vector(-1,-1){2}}
\put(10.,10.){\vector(1,-1){2}}
\put(3.5,6.){\vector(0,-1){0.5}}
\put(12.,6.){\vector(0,-1){0.5}}
\put(3.5,4.5){\vector(0,-1){1.}}
\put(12.,4.5){\vector(0,-1){1.}}
\put(3.750,3.75){\framebox(1.75,0.5){$j=3$}}
\put(12.25,3.75){\framebox(1.75,0.5){$j=3$}}
\end{picture}

\renewcommand{\arraystretch}{1.0}
  \caption{An example of
  $\pi$\hyph representative determination
  by Algorithm~\ref{pipredst}.}\label{ppipredst}
\end{figure}

The Algorithm~\ref{pipredst} is not efficient for extremely symmetric
matrices, such as
$I_n$:
in that case bound step does not ever occur,
because all the remaining
rows are always equally good.
Hence, Algorithm~\ref{pipredst} must be improved,
in order to detect some
symmetries, and to avoid some unnecessary repetitions.
Suppose that there remain $l$ rows not included in $A_\pi$,
 and that the column classes defined by $\Sigma^{(n-l-1)}$
 are such, that
all column classes in the remaining rows are \emph{uniform}
(they contain either all ones or all zeros),
except for at most one column class, which in that case has
$l$ columns,
with the row and column sums both equal to
$l-1$ or $1$.
Then, because of the symmetry,
it is enough
to put in $L^{(n-l-1)}$ only one of the $l$ remaining rows.
After the incorporation of this simple heuristic, the algorithm
much more efficiently deals with the
matrices such as
$I_n$, the complement of $I_n$,
and the other highly symmetric matrices.

Using Algorithm~\ref{pipredst}, it is possible to determine
$A_\phi$ for
an arbitrary $A\in\mathcal{A}_n$:
it is enough to find $\pi$\hyph representatives of all
$(n+1)^2$ matrices
$X_i A X_j$, $0\le i,j\le n$,
and then to choose the smallest among them.

One of the outputs from Algorithm~\ref{pipredst}
is the number of the
pairs of row/column permutations,
transforming $A$ into $A_\pi$.
That number is used to determine the size of the $\pi$\hyph class
of  $A^T$,
as it will be demonstrated below.

Consider the problem of counting the
matrices in the $\pi$\hyph class of an arbitrary
$A\in\mathcal{A}_n$.
For an arbitrary
$B\in\mathcal{A}_n$ let $B_0$ denote the matrix obtained from $B$
by sorting its rows.
If $A$ has $i_k$ groups of $k$ equal rows, $1\le k\le n$, then
the number of matrices that could be obtained from $A$
by row permutations is
\[a=n!/\prod_{k=1}^n {i_k}!\]
The representative of these $a$ matrices is $A_0$.
An arbitrary matrix $A'$, obtained from $A$ by a column permutation,
generates in the same manner a new set of $a$ matrices if and only if $A'_0\neq A_0$.
If the number of different matrices $A_0'$
is $b$, then the size of the
$\pi$\hyph class of $A$ is $ab$.
It is simpler to obtain $b$
by counting the number $p$ of column permutations $A'$ of $A$
satisfying $A_0'=A_0$, because
$b=n!/p$.
Note that $p$ is preserved by row and column permutations of $A$.

Applying
Algorithm~\ref{pipredst} to $(A^T)_\pi$,
$p$ is obtained even more easily.
Indeed, suppose that $A$ is already a $\pi$\hyph representative, i.e.
$A=A_\pi$. Then Algorithm~\ref{pipredst}
counts the row permutations $A''$ of $A$,
such that there exists a column permutation $A'''$ of $A''$, equal to
$A_\pi$.
Now we find $A'=((A^T)_\pi)^T$ and apply Algorithm~\ref{pipredst}
(again) to $(A')^T$.
The matrix $(A')^T$ is a $\pi$\hyph representative, because
$((A')^T)_\pi=(A')^T$.
Algorithm~\ref{pipredst} gives
the number of
row permutations $(A'')^T$ of $(A')^T$,
such that there exists a column permutation $(A''')^T$ of $(A'')^T$,
equal to
$(A')^T$.
In other words, we obtain the number of
column permutations $A''$ of $A'$,
such that there exists a row permutation $A'''$ of $A''$, equal to
$A'$ --- which is exactly $p$ ($count$ in Algorithm~\ref{pipredst}).

\begin{exmp}
Looking again at Example~\ref{ex2},
we see that there are two pairs
$(P,Q)$ that minimize $PAQ$.
Therefore, there are $3!^2/2=18$ matrices in the
$\pi$\hyph class of $A^T$.
\end{exmp}

The problem of counting the
matrices in the $\SNF$\hyph class of an arbitrary
$A\in\mathcal{A}_n$ is much harder.
It is even harder is to enumerate the sets
$\mathcal{A}_{n,k}=\{A\in\mathcal{A} \mid \rang A=k\}$, $0\le k\le n$:
(especially  $m_n=\mathcal{A}_{n,n}$)
 We now explicitly enumerate the sets
 $\mathcal{A}_{n,1}$, $\mathcal{A}_{n,2}$,
using the following characterization of
matrices in $\mathcal{A}_{n,2}$.

\begin{lem}\label{rang2}
If the matrix
$A\in\mathcal{A}_{n,2}$ contains three different
nonzero columns
$a$, $b$, $c$,
then one of them is equal to the sum of the other two,
for example
$c=a+b$.
Furthermore, the set of nonzero rows of the matrix $[a\ b]$ equals to
$\{[0\ 1],[1\ 0]\}$.
There can not be four different nonzero columns in $A$.
\end{lem}

\begin{pf}
Suppose $A\in\mathcal{A}_{n,2}$.
If two nonzero columns of $A$ are linearly dependent,
then they are obviously equal.
Suppose $a$, $b$, $c$ are the three different nonzero
linearly dependent columns, i.e.
$\alpha a+\beta b+\gamma c=0$
for some integers
$\alpha$, $\beta$, $\gamma$.
The coefficients
$\alpha$, $\beta$, $\gamma$
must be nonzero; otherwise, if for example $\alpha=0$,
then
$\beta b+\gamma c=0$ implies $b=c$.
Denote by $U$ the set
of nonzero rows of the $n\times 3$ matrix
$[a\ b\ c]$.
Then
\begin{itemize}
\item
$|U|>1$; otherwise it would be $a=b=c$.
\item
$U\cap\{[1\ 0\ 0],[0\ 1\ 0],[0\ 0\ 1]\}=\emptyset$;
if, for example $[1\ 0\ 0]\in U$,
then $\alpha=0$.
\item
therefore, $U\subseteq \{[1\ 1\ 1],[0\ 1\ 1],[1\ 0\ 1],[1\ 1\ 0]\}$ and
$U\neq \{[1\ 1\ 1]\}$.
\item
$[1\ 1\ 1]\notin U$;
if
$[1\ 1\ 1]\in U$, and for example
$[0\ 1\ 1]\in U$,
then from
$\alpha+\beta+\gamma=0$ and
$\beta+\gamma=0$, it follows
$\alpha =0$.
\item
$U\neq \{[0\ 1\ 1],[1\ 0\ 1],[1\ 1\ 0]\}$;
otherwise
$\beta+\gamma=0$,
$\alpha+\gamma=0$,
$\alpha+\beta=0$
implies
$\alpha=\beta=\gamma=0$.
\end{itemize}
Hence, there are three possibilities for $U$ left:
$\{[0\ 1\ 1],[1\ 0\ 1]\}$, or
$\{[0\ 1\ 1],[1\ 1\ 0]\}$, or
$\{[1\ 0\ 1],[1\ 1\ 0]\}$.
If
$U=\{[0\ 1\ 1],[1\ 0\ 1]\}$, then
$\beta+\gamma=0$,
$\alpha+\gamma=0$
implies
$(\alpha,\beta,\gamma)=\gamma(-1,-1,1)$,
i.e. $c=a+b$;
the set of nonzero rows of $[a\ b]$ is
$\{[0\ 1],[1\ 0]\}$.
The two other cases are symmetrical.

Suppose that $A$ contains four different columns
$a$,
$b$,
$c$,
$d$.
Then we must have for example $c=a+b$ and
the set of nonzero rows of $[a\ b]$ is
 $\{[0\ 1],[1\ 0]\}$.
Applying the first part of Lemma to
$a$, $b$, $d$,
we conclude that
$d=a+b$ or
$a=b+d$ or
$b=a+d$.
But
$d=a-b$ and
$d=b+a$ are impossible, 
and
$d=a+b$ implies $d=c$.
The lemma is proved.
\qed \end{pf}

\begin{thm}
a) For an arbitrary
$A\in\mathcal{A}_n$ the following three statements are equivalent:
\begin{enumerate}
\item
$\rang A=1$;
\item
$\SNF(A)=(1,0^{n-1})$;
\item
$A$ contains a column $a\neq 0$,
such that all nonzero columns of $A$ are equal to $a$.
\end{enumerate}
The number of matrices in $\mathcal{A}_{n,1}$ equals
\[ |\mathcal{A}_{n,1}|=(2^n-1)^2.\]

b) For an arbitrary
$A\in\mathcal{A}_n$ the following three statements are equivalent:
\begin{enumerate}
\item
$\rang A=2$;
\item
$\SNF(A)=(1,1,0^{n-2})$;
\item
\begin{itemize}
\item
$A$ contains the two nonzero columns $a\neq b$,
such that all columns of $A$ are in $\{0,a,b\}$,
or
\item
$A$ contains the two nonzero columns $a\neq b$,
such that the set of
nonzero rows of $[a\ b]$ equals
$\{[0\ 1],[1\ 0]\}$, and that the set of
nonzero columns of $A$ is $\{a,b,a+b\}$.
\end{itemize}
\end{enumerate}
The number of matrices in $\mathcal{A}_{n,2}$ equals
\[ |\mathcal{A}_{n,2}|=(3^n-2\cdot2^n+1)(2\cdot4^n-3\cdot3^n+1)/2.\]
\end{thm}

\begin{pf}
a) If $\rang A=1$
then $A$ contains
nonzero column $a$, such that all nonzero columns of $A$ are
equal to $a$.
By subtracting one of nonzero columns from the others,
we obtain an
equivalent matrix with exactly one nonzero column $a$.
By the column permutation column $a$ is moved to the first position,
and by the row permutation some $1$ is moved to the upper left corner.
By subtracting the first row from the other nonzero rows,
we obtain that
SNF of $A$ is $(1,0^{n-1})$.
How many matrices of rank $1$ there are?
The number of choices for nonzero column $a$ is
$2^n-1$, and the number of
matrices corresponding to the fixed $a$ is
$2^n-1$: each its column is $0$ or $a$, but at least one of them
has to be equal to $a$. Hence,
$|\mathcal{A}_{n,1}|=(2^n-1)^2.$

b) If $\rang A=2$ then $A$ contains two linearly independent
columns,
such that the other
columns are their linear combinations. %
The number of different nonzero
columns in $A$ is either two or it is greater than two.

\begin{description}
\item[Case 1.]
Suppose there are exactly two different nonzero columns $a$, $b$ in $A$.
The number of such matrices $A$ is
\[
\binom{2^n-1}{2}
(3^n-2\cdot 2^n+1).
\]
Indeed, the number of choices for $a$, $b$ equals to the
above binomial coefficient.
Without loss of generality we suppose that $a<b$.
For fixed $a$, $b$, by the inclusion-exclusion principle
the number of matrices $A$ is
$3^n-2\cdot 2^n+1$,
because
\begin{itemize}
\item
$3^n$ is the number of matrices with the columns from the set
$\{0,a,b\}$,
\item
$2^n$ is the number of matrices without $a$, and also
the number of matrices without $b$,
\item
$1$ is the number of matrices without $a$ and $b$.
\end{itemize}
\item[Case 2.] If there are more than two
different nonzero columns in $A$, then
by Lemma~\ref{rang2}
there are
two different nonzero columns $a$, $b$  ($a<b$)
in $A$,
such that the set of
nonzero columns in $A$ is
$\{a,b,c=a+b\}$,
and such that the row set of the matrix $[a\ b]$ is
$\{[0\ 1],[1\ 0]\}$.
There are
$(3^n-2\cdot 2^n+1)/2$
choices for columns $a$, $b$ satisfying these conditions.
Indeed,
consider all matrices $[a\ b]$, $[b\ a]$:
\begin{itemize}
\item
$3^n$ is the number of matrices with the row set
$\{[0\ 0],[0\ 1],[1\ 0]\}$,
\item
$2^n$
is the number of matrices without the row
$[0\ 1]$, and also the number of matrices without the row
$[1\ 0]$,
\item
$1$
is the number of matrices without the rows
$[0\ 1]$,
$[1\ 0]$).
\end{itemize}
The number of matrices $[a\ b]$ is therefore
$(3^n-2\cdot 2^n+1)/2$.
The number
$4^n-3\cdot 3^n+3\cdot 2^n-1$
of matrices with the set of nonzero columns
$\{a,b,c\}$ (where $c=a+b$) is also obtained by the
inclusion-exclusion principle:
\begin{itemize}
\item
$4^n$ is the number of matrices with all the columns $0$, $a$, $b$, $c$;
\item
$3^n$ is the number of matrices without the column $a$
(and analogously without $b$, $c$);
\item
$2^n$ is the number of matrices without columns $a$, $b$
(and analogously without $a$, $c$; and without $b$, $c$);
\item
$1$ is the number of matrices without columns $a$, $b$, $c$.
\end{itemize}
Therefore, the number of matrices of the rank $2$,
with more than two different nonzero columns equals
\[
(3^n-2\cdot 2^n+1)(4^n-3\cdot 3^n+3\cdot 2^n-1)/2.
\]
The total number of matrices in $\mathcal{A}_{n,2}$ equals
\begin{eqnarray*}
& &\left(2\binom{2^n-1}{2}+(4^n-3\cdot 3^n+3\cdot 2^n-1)\right)
(3^n-2\cdot 2^n+1)/2=\\
&=&(3^n-2\cdot2^n+1)(2\cdot4^n-3\cdot3^n+1)/2.
\end{eqnarray*}
\end{description}

In either case, in order to obtain $\SNF(A)$,
the other nonzero columns are first transformed to $0$ by subtracting
$a$, $b$ or $a+b$
from them.
Next, in $[a\ b]$ there is a
row $[0\ 1]$, because $a<b$;
using that $1$, the other elements of $b$ are changed to $0$.
Finally, choosing some $1$ in $a$,
and subtracting if necessary that row from the others,
after permuting rows/columns, we obtain the SNF.
Hence, $\rang A=2$ implies
$\SNF(A)=(1,1,0^{n-2})$.
\qed \end{pf}

\subsection{Iterative classification of \nj matrices}

According to Lemma~\ref{pobrub} we have 
\[
\mathcal{A}_{n+1}^\pi = \cup_{A\in\mathcal{A}_n^\pi} \bord_\pi(A).
\]
By
changing the order of calculations, it is possible to simplify
repeated determination of $\pi$\hyph representatives of matrices
from $\bord(A)$ by~Algorithm~\ref{pipredst}.
Matrices $B$ in $\bord(A)$ are of the
form~\eqref{bordA}.
For each $y$ the  $\pi$\hyph representatives of $B$'s
corresponding to various inserted rows $[x\ b]$ are found
spending smaller number of steps.
The point is that the rows of the $\pi$\hyph representative preceding
the row $[x\ b]$ are already determined for some previous
variants for that row.

Somewhat more detailed description follows.
Determine first the
$\pi$\hyph representative of the matrix, corresponding to $x=0$,
$b=0$; the inserted zero row $[x\ b]$ is certainly the first row
in the $\pi$\hyph representative. The corresponding row and column
permutations $P$, $Q$ are recorded.
The remaining pairs $(x,b)$ are then considered in turn,
lexicographically ordered.
The question arises, to which position
$l$ might $[x\ b]$ be moved during the
$\pi$\hyph representative determination,
skipping the determination of first $l-1$ rows of the representative.
The obvious lower bound for $l$
is the smallest among all positions where
the previous rows $w$,
obtained from $[x\ b]$ by changing
exactly one $1$ into $0$, have been moved
(except if there was an alternative to
$w$ during that step,
i.e. if $L^{(i)}$ had more than one member at the moment
when $w$ arrived to its destination).


Instead of extending all $A\in\mathcal{A}_n^\pi$,
it is enough to extend 
the matrices from the set
$\mathcal{A}_n^\phi$ of all $\phi$\hyph representatives in \ann.
By extending all $A\in\mathcal{A}_n^\phi$
a subset of $\mathcal{A}_{n+1}^\pi$ is obtained;
the set of
$\phi$\hyph representatives of matrices from that subset
is exactly $\mathcal{A}_{n+1}^\phi$.

It is convenient to use a balanced tree
to collect
$\pi$\hyph representatives in an ordered fashion.
We chose AVL tree~\cite{AVL} --- the binary search tree
satisfying the condition that, for every node,
the difference between the heights of its left and right subtrees
is at most $1$.
For $n=8$, in order to save memory,
a combination of AVL tree and the sorted array of matrices is used:
from time to time the content of the tree is merged into the array.
After collecting all $\pi$\hyph representatives,
the $\pi$\hyph representatives set is reduced to the
corresponding $\phi$\hyph representatives set.
To determine the set of $\phi$\hyph representatives, 
corresponding to
a given set of
$\pi$\hyph representatives, the following simple algorithm
is used.

\begin{alg}\label{phipredst}
Reduction of a given set $L_\pi$ of $\pi$\hyph representatives
  to the set $L_\phi$ of corresponding
  $\phi$\hyph representatives.
\begin{tabbing}
 \quad \= \quad \= \quad   \kill
 \{ $T$ --- auxiliary AVL tree used to collect $\pi$\hyph representatives. \} \\
 \WHILE $L_\pi\neq \emptyset$ \\
 \> \WHILE there is a space in $T$ for at least $(n+1)^2$ matrices\\
 \> \>  remove the first matrix $A$ from $L_\pi$;\\
 \> \>  generate the set $T_A$ of $\pi$\hyph representatives contained in the $\phi$\hyph class of $A$;\\
 \> \>  insert $T_A$ into $T$;\\
 \> \>  insert $A_\phi$ into $L_\phi$;\\
 \>  remove from $L_\pi$ all the matrices contained in $T$;\\
 \>  $T\la \emptyset$;\\
\end{tabbing}
\end{alg}

The classification of \ao
lasted about a month in parallel on five
PC's.
A huge number of collected
$\pi$\hyph representatives of order
$n=8$ caused serious difficulties.
The space requirement is reduced by dividing $\pi$\hyph representatives into
subsets, according to their SNF.
For each extended matrix, its SNF is determined,
and the
$\pi$\hyph representatives are classified into subsets
with the same SNF.
These subsets are then independently
processed.
The hardest was the SNF-class
$(1^7,0)$, with
$5204144555$ $\pi$\hyph representatives contained in
a number of non disjoint subsets.
These subsets were independently processed
by Algorithm~\ref{phipredst}, producing the
non disjoint sets of
$\phi$\hyph representatives; their union consists of
$71348129$ $\phi$\hyph representatives,
approximately $1/3$ of matrices in $\mathcal{A}_8^\phi$.

In order to save the space,
$L_\pi$ and $L_\phi$ are stored in a sorted, compressed form:
one byte for each matrix row;
the group of consecutive
matrices with the same first $n-2$ rows
is stored so that the common
$n-2$ rows are stored only once.
As a result, the average space for a matrix
of order $8$ was little more than two bytes.

If somebody tries to extend $\phi$\hyph representatives of order $8$,
he could expect to process about
$300$ times more $\phi$\hyph representatives,
each giving approximately
 $4$ times more $\pi$\hyph representatives.
Therefore, the classification of matrices
 of order $9$ is expected to last $1000$ times longer,
requiring huge memory.

\subsection{Results of classification}\label{rez}

We start with the simplest nontrivial case.

\begin{exmp}
The $16$ matrices of order $2$ are divided into $3$ $\phi$\hyph classes,
which are further subdivided into $7$ $\pi$\hyph classes:
{\scriptsize
\begin{eqnarray*}
& &\left\{
 \left\{
 \left[\begin{array}{rr}  0 & 0 \\  0 & 0\end{array}\right]
 \right\}
\right\}, \\
& &\left\{
 \left\{
 \left[\begin{array}{rr}  0 & 0 \\  0 & 1\end{array}\right],
 \left[\begin{array}{rr}  0 & 0 \\  1 & 0\end{array}\right],
 \left[\begin{array}{rr}  0 & 1 \\  0 & 0\end{array}\right],
 \left[\begin{array}{rr}  1 & 0 \\  0 & 0\end{array}\right]
\right\},
\left\{
 \left[\begin{array}{rr}  0 & 0 \\  1 & 1\end{array}\right],
 \left[\begin{array}{rr}  1 & 1 \\  0 & 0\end{array}\right]
\right\},\right.\\
& &\quad \left.\left\{
 \left[\begin{array}{rr}  0 & 1 \\  0 & 1\end{array}\right],
 \left[\begin{array}{rr}  1 & 0 \\  1 & 0\end{array}\right]
\right\},
 \left\{
 \left[\begin{array}{rr}  1 & 1 \\  1 & 1\end{array}\right]
 \right\}
\right\}, \\
& &\left\{
 \left\{
 \left[\begin{array}{rr}  0 & 1 \\  1 & 0\end{array}\right],
 \left[\begin{array}{rr}  1 & 0 \\  0 & 1\end{array}\right]
\right\},
\left\{
 \left[\begin{array}{rr}  0 & 1 \\  1 & 1\end{array}\right],
 \left[\begin{array}{rr}  1 & 0 \\  1 & 1\end{array}\right]
 \left[\begin{array}{rr}  1 & 1 \\  0 & 1\end{array}\right],
 \left[\begin{array}{rr}  1 & 1 \\  1 & 0\end{array}\right]
\right\}
\right\}.
\end{eqnarray*}
}
\end{exmp}

In Table~\ref{pi3} all the $36$ $\pi$\hyph representatives
of order $3$ are shown.
The $5$ SNF-classes are in separate  blocks, divided into
compartments with $\phi$\hyph classes. The first matrix in each
$\phi$\hyph class is the smallest $\pi$\hyph representative,
i.e. the $\phi$\hyph representative.
For each
$\pi$\hyph  and SNF-class, their size is given.
The matrices are represented by hexadecimal vectors,
each component representing a row of a matrix.
For example, the last
vector $(3,5,6)$ in Table~\ref{pi3} represents the matrix
{\small
\[
\left[\begin{array}{rrr}
  0 & 1 & 1 \\
  1 & 0 & 1 \\
  1 & 1 & 0
\end{array}\right].
\]
}
The matrix $(1,2,5)$ is a
$\pi$\hyph representative of the matrix from
Example~\ref{ex2}.

In Table~\ref{phi4} all the $39$ $\phi$\hyph representatives
of order $4$
are shown, together
with the sizes of their
$\phi$-classes.

In Table~\ref{brojanje}
$\rho_n$,
$|\mathcal{A}_n^\phi|$,
$s_n$, $a_n$,
$|\mathcal{D}_n|$,
 and the set $\mathcal{D}_n$ are given
for $1\le n\le 8$, where
$s_n=|\mathcal{S}_n|$ and
$\rho_n=(2^{n^2}/(n+1)!^2)/|\mathcal{A}_n^\phi|$.
In the last row of Table~\ref{brojanje}
 $s_9$,
$|\mathcal{D}_9|$, $a_9$,
 $\mathcal{D}_9$ are given;
the explanation how they are obtained will be given in section~\ref{snf9}.

\begin{table}
 \caption{The numbers of equivalence classes in \ann.}
 \label{brojanje}
{\small
\[
\begin{array}{||r|r|r|r|r|r|l||}
\hline\hline
n
& \rho_n
&|\mathcal{A}_n^\phi|&s_n& a_n &|\mathcal{D}_n| &\multicolumn{1}{c||}{\mathcal{D}_n} \\
\hline\hline
1 &0.250&        2&    2&     2&   2&\left\{0,1\right\}                       \\\hline
2 &0.148&        3&    3&     2&   2&\left\{0,1\right\}                       \\\hline
3 &0.074&       12&    5&     3&   3&\left\{0-2\right\}                       \\\hline
4 &0.117&       39&    8&     4&   4&\left\{0-3\right\}                       \\\hline
5 &0.167&      388&   14&     6&   6&\left\{0-5\right\}                       \\\hline
6 &0.334&     8102&   26&    10&  10&\left\{0-9\right\}                       \\\hline
7 &0.528&   656103&   56&    19&  22&\left\{0-18,20,24,32\right\}             \\\hline
8 &0.701&199727714&  129&    41&  46&\left\{0-40,42,44,45,48,56\right\}       \\\hline
9 &     &         &  333&   103& 114&\left\{0-102,104,105,108,110,\right.\\
  &     &         &     &      &    &\left. 112,116,117,120,125,128,144\right\}   \\
\hline\hline
\end{array}
\]
}\end{table}

Denote by $F(n)$ the following statement:
\begin{eqnarray}\label{Fn}
 & &%
\text{$A\in\mathcal{A}_n$, satisfying
$\SNF(A)=d=(d_1,d_2,\ldots,d_n)$
exists if and only if}\\
& &\text{there exists $A'\in\mathcal{A}_{n+1}$, satisfying
$\SNF(A')=d'=(d_1,d_2,\ldots,d_n,0)$.}
\nonumber
\end{eqnarray}
Obviously, the first condition implies the second one.
The implication in the opposite direction is not obvious at all;
it would follow from the following
stronger statement:
\begin{eqnarray*}
 & &H(n):
\text{ Let
$A'\in\mathcal{A}_{n+1}$, $\rang A'=n$,
and $\SNF(A')=d'=(d_1,d_2,\ldots,d_n,0)$.}\\
& &\text{Then $A'$ has at least one minor $A\in\mathcal{A}_n$
with $\SNF(A)=d=(d_1,d_2,\ldots,d_n)$.}
\end{eqnarray*}
But the following matrix $F\in\mathcal{A}_{10}$ is a
counterexample to $H(10)$:
\[
F=
\left[\begin{array}{rr}
A&B\\
C&D
\end{array}\right]=
\left[\begin{array}{rrrr@{\,}|@{\,}rrrrrr}
0&0&1&1&1&0&1&0&0&1 \\
0&1&1&0&0&1&0&1&0&1 \\
1&1&0&0&0&1&1&0&1&0 \\
1&0&0&1&1&0&0&1&1&0 \\
\hline
0&0&1&1&0&0&0&1&1&1 \\
1&1&0&0&0&0&1&1&1&0 \\
1&0&1&0&0&1&1&1&0&0 \\
0&1&0&1&1&1&1&0&0&0 \\
0&1&1&0&1&1&0&0&0&1 \\
1&0&0&1&1&0&0&0&1&1
\end{array}\right].
\]
The matrix $F$ consists of blocks $A,B,C,D$,
having $2,3,2,3$ ones in each row respectively,
and also having
$2,2,3,3$ ones in each column, respectively;
$F$ is singular, because the sums of rows of
$[A\ B]$ and $[C\ D]$ are equal.
It can be verified that
$\rang F=9$,
$\SNF(F)=(1^9,0)$,
but all minors of $F$ have SNF different from
$(1^9)$.

In Table~\ref{kat8} the SNF-representatives of
matrices in $\mathcal{A}_n$,
$n\le 8$, are listed,
accompanied with the size measures of corresponding SNF-classes
(the number of matrices, the number of
$\pi$\hyph representatives and the number of
$\phi$\hyph representatives in each SNF-class).
The sizes of $\pi$\hyph classes
are determined
using Algorithm~\ref{pipredst}.
The classes are ordered lexicographically
by the SNF (with zeros moved to the end of SNF).

One can verify this classification
starting from the sorted list of all $\phi$\hyph representatives.
For each of them one has to check if it is indeed a
$\phi$\hyph representative. The next step is to sum the numbers of
$\pi$\hyph representatives in all $\phi$\hyph classes,
and to compare the sum with the corresponding entry in
Table~\ref{ubas}.
One could also check that the sum of sizes of SNF-classes in \an
equals $2^{n^2}$ for each $n\le 8$,
see Table~\ref{kat8}.
The sorted lists of $\phi$\hyph representatives for $n\le 8$
can be downloaded from
http://www.matf.bg.ac.yu/~ezivkovm/01matrices.htm.

We now review some interesting facts, 
which are seen from Table~\ref{kat8}.

Let $T(n,k)=|\mathcal{A}_{n,k}|$.
In Table~\ref{apsrang} the numbers $T(n,k)$,  $0\le k\le
n\le 8$, are shown
(of course, they are easily obtained from Table~\ref{kat8}). 
The part of Table~\ref{apsrang} corresponding
to $n\le 7$ is the same as in~\cite{Ziegler}; it is also an entry
in~\cite[Sequence A064230]{Slo}. Another interesting entry
in~\cite[Sequence A055165]{Slo} is the sequence $m_n$, where $m_n$
is the number of regular \nj matrices of order $n$ --- the
diagonal of Table~\ref{apsrang}. 
The seemingly new member of that sequence
is $m_8=10160459763342013440$. If we suppose that all matrices in
\an are equiprobable, then the rank probability distribution is
shown in Table~\ref{verrang} for $n\le 8$. Looking at
Table~\ref{verrang}, one could erroneously conclude that
large fraction of matrices in \an is singular. In fact, the
fraction of singular matrices in \an tends to $0$ 
for $n$ large~\cite{Komlos}.

\begin{table}
\caption{The number of matrices of the rank $k$ in \ann, $n\le 8$.}
\label{apsrang}
{\small

\begin{tabular}{||r||r@{\ry}r@{\ry}r@{\ry}r@{\ry}r@{\ry}r@{\ry}r@{\ry}r||}
\hline\hline
~~n  & 1&2&3  & 4   &    5    &        6   &        7       &             8      \\
k~~  &  & &   &     &         &            &                &                    \\
\hline\hline
 0   & 1&1&1  &1    &1        &1           & 1              &                   1\\
 1   & 1&9&49 &225  &961      &3969        & 16129          &               65025\\
 2   &  &6&288&6750 &118800   &1807806     & 25316928       &           336954750\\
 3   &  & &174&36000&3159750  &190071000   & 9271660734     &        397046059200\\
 4   &  & &   &22560&17760600 &5295204600  & 1001080231200  &     144998212423680\\
 5   &  & &   &     &12514320 &34395777360 & 32307576315840 &   17952208799918400\\
 6   &  & &   &     &         &28836612000 & 259286329895040&  720988662376725120\\
 7   &  & &   &     &         &            & 270345669985440& 7547198043595392000\\
 8   &  & &   &     &         &            &                &10160459763342013440\\
\hline\hline
\end{tabular}
}
\end{table}

\begin{table}
\caption{The probability that a random matrix in \an has the rank
$k$, $0\le k\le n\le 8$.} \label{verrang}
{\small
\begin{tabular}{||r||rrrrrrrr||}
\hline\hline
~~n&  1&  2     &  3      &   4     &    5    &     6   &  7      &   8     \\
k~~&   &        &         &         &         &         &         &         \\
\hline\hline
 0&0.5 & 0.0625 & 0.00195 & 0.00002 & 0.00000 & 0.00000 & 0.00000 &0.00000  \\
 1&0.5 & 0.5625 & 0.09570 & 0.00343 & 0.00003 & 0.00000 & 0.00000 &0.00000  \\
 2&    & 0.3750 & 0.56250 & 0.10300 & 0.00354 & 0.00003 & 0.00000 &0.00000  \\
 3&    &        & 0.33984 & 0.54932 & 0.09417 & 0.00277 & 0.00002 &0.00000  \\
 4&    &        &         & 0.34424 & 0.52931 & 0.07706 & 0.00178 &0.00001  \\
 5&    &        &         &         & 0.37296 & 0.50052 & 0.05739 &0.00097  \\
 6&    &        &         &         &         & 0.41963 & 0.46059 &0.03908  \\
 7&    &        &         &         &         &         & 0.48023 &0.40913  \\
 8&    &        &         &         &         &         &         &0.55080  \\
\hline\hline
\end{tabular}
}
\end{table}

It turns out that
$F(n)$~\eqref{Fn} is true for $n\le7$, i.e.
the set of SNF's of rank $k$ is the same for all $n$,
$k\le n\le 8$.
For example,
the SNF-representative of the SNF-class
$(1,1,2,0^{n-3})$ is the matrix
$(0^{n-3},3,5,6)$ for $3\le n\le 8$.

The smallest $n$ for which there are two
matrices in \an with the same determinant,
but with different
SNF's is $5$:
$\SNF(\text{3,C,15,16,19})=(1,1,1,4)$
and $\SNF(\text{3,5,9,11,1E})=(1, 1, 2, 2)$.

In Table~\ref{velpxkl} the possible numbers of
$\pi$\hyph orbits inside
$\phi$\hyph orbits are shown
for $1\le n\le 8$.
These numbers are between $1$ and $(n+1)^2$;
as it is seen,
the value $(n+1)^2$ is attained only if $n\ge 5$.

\begin{table}
 \caption{The possible numbers of $\pi$\hyph orbits inside $\phi$\hyph orbits
of \ann.}
 \label{velpxkl}
{\small
\[
\begin{array}{||r||l||}
\hline\hline
n & \text{The set of $\phi$\hyph orbit sizes} \\
\hline\hline
1& \{1 \} \\
2& \{1,2,4 \} \\
3& \{1,2,4,5,9 \} \\
4& \{1-5,7,9-11,13,16,17 \} \\
5& \{1-18,20,21,25,26,30,36 \} \\
6& \{1,2,4-27,29-32,35-37,42,49 \} \\
7& \{1-38,40,42-44,48-50,56,64 \} \\
8& \{1-46,48-51,53,54,56-58,63-65,72,81 \} \\
\hline\hline
\end{array}
\]
}\end{table}

If $A\in\mathcal{A}_n$, $A\sim I_n$ and $B\in\bord(A)$,
then $\SNF(B)$
contains at least $n$ ones, see \eqref{zz}.
The question arises,
what are the possible values of the last element of $\SNF(B)$,
i.e. which values can take $|\det B|$?
The largest possible values of $|\det B|$
under these assumptions, along with the examples of
matrices $B$ for which these values are attained,
are given in Table~\ref{extremi}.
In fact, the matrices from Table~\ref{extremi} maximize 
 $|\det B/\det A|$ for all regular
$A\in\mathcal{A}$, $n\le 8$.

\begin{table}
 \caption{The maximal ADV's of matrices from $\mathcal{A}_{n+1}$,
 obtained by extending matrices equivalent to $I_n$.}
 \label{extremi}
{\small
\[
\begin{array}{||r|r||rrrrrrrrr||}
\hline\hline
n & |\det A| & \multicolumn{9}{c||}{A}\\
\hline\hline
3 & 3 &  3 &  5 &  9 &  E &    &    &     &     &      \\
4 & 5 &  3 &  5 &  E & 16 & 19 &    &     &     &      \\
5 & 9 &  3 &  D & 15 & 1A & 26 & 39 &     &     &      \\
6 & 18&  7 & 19 & 2A & 34 & 4C & 53 & 65  &     &      \\
7 & 40&  7 & 19 & 2A & 56 & 65 & 9C & B3  & CB  &      \\
8 &105&  7 & 39 & 5A & AC & D5 & E3 &136  &14D  & 19B  \\
\hline\hline
\end{array}
\]
}\end{table}

More generally, it is interesting to describe the relationship
of
$\SNF(A)$ to $\SNF(A')$ if
$A'\in \bord(A)$.
During iterative classification,
the sets
\[
\{\SNF(B) \mid B\in \bord(A),\ A\in \mathcal{A}_n,\ \SNF(A)=s\}
\]
are recorded for all SNF-classes $s\in\mathcal{S}_n$.
The results are represented by the incidence matrix $M_n$
of dimensions
$|\mathcal{S}_{n}|\times|\mathcal{S}_{n+1}|$,
with entries
\begin{equation}\label{SNFinc}
m_{s,s'}=
\left\{
 \begin{array}{ll}
  $1$, & \text{ if there exist $A\in \mathcal{A}_n$ and $B\in \bord(A)$,
  with SNF's $s$ and $s'$ }\\ 
  $0$, & \text{ otherwise}
 \end{array}
\right. .
\end{equation}

Let $G(n)$, denote  the following statement:
\begin{eqnarray}\label{Gn}
 & & G(n):
\text{ There exist matrices $A\in\mathcal{A}_n$,
$A'\in\bord(A)$, such that } \\
& &\text{$\SNF(A)=(d_1,d_2,\ldots,d_n)$,
$\SNF(A')=(d'_1,d'_2,\ldots,d'_n,d'_{n+1})$
if and only if }\nonumber\\
& &\text{there exist matrices $B\in\mathcal{A}_{n+1}$,
$B'\in\bord(B)$, such that}\nonumber \\
& &\text{ $\SNF(B)=(d_1,d_2,\ldots,d_n,0)$,
$\SNF(B')=(d'_1,d'_2,\ldots,d'_n,d'_{n+1},0)$.}
\nonumber
\end{eqnarray}
By exhaustive search it is verified that $G(n)$ is true
for $n\le 6$, enabling
to put all the transposed incidence
matrices $M_n$,
$n\le 7$ together
into single Table~\ref{SNFincid}.
The $1$'s are represented by $\bullet\,$;
the $0$'s are represented by $\star$
if they are the consequence of the following Lemma
(describing constraints for $\SNF(A')$ if $A'\in\bord(A)$);
otherwise, they are represented by
$\circ\,$.

\begin{lem}\label{SNFogr}
For an arbitrary $A\in\mathcal{A}_n$, let
$A'\in\bord(A)$, and let
$\SNF(A)=(d_1,d_2,\ldots,d_n)$,
$\SNF(A')=(d'_1,d'_2,\ldots,d'_n,d'_{n+1})$.
Then
\begin{enumerate}
  \item $\rang A\le \rang A'\le \rang A +2$;
  \item $d'_1d'_2\ldots d'_i$ divides $d_1d_2\ldots d_i$ for all $i$,
$1\le i\le \rang A$;
  \item $\prod_{i=1}^{n-1}d_i$ divides
$\det A'$.
\end{enumerate}
\end{lem}

\begin{pf}

\begin{enumerate}
  \item
The first inequality follows from the fact that
the rank of a submatrix is a lower bound on the rank
of a matrix.  The second inequality follows from
the observation that $A'$ is an at most rank 2 perturbation
of $A$.
  \item This is a direct consequence of the fact that
  $d_1'd_2'\ldots d_i'$
is the largest common divisor of all minors of $A'$ of order $i$,
see for example~\cite{gantmaher}.
  \item
  Let $P$, $Q$ be the matrices such that
$\SNF(A)=PAQ=D=(d_1,d_2,\ldots,d_n)$,
$|\det P|=|\det Q|=1$.
Let
\[A'=
\left[\begin{array}{rr}
  A & y \\
  x & b
\end{array}\right].
\]
The case $\det A'=0$ is trivial; suppose $\det A'\neq 0$.
If
$xQ=[a_1\ a_2\ \ldots a_n]$,
$Py=[c_1\ c_2\ \ldots c_n]^T$,
then from the identity
\[
\left[\begin{array}{rr}
  P & 0 \\
  0 & 1
\end{array}\right]
\left[\begin{array}{rr}
  A & y \\
  x & b
\end{array}\right]
\left[\begin{array}{rr}
  Q & 0 \\
  0 & 1
\end{array}\right]=
\left[\begin{array}{rr}
  D & Py \\
  xQ & b
\end{array}\right]
\]
it follows 
(another way to express determinants 
of matrices obtained by extension, see~\eqref{detobr})
\begin{equation}\label{z}
\det A'=
b\prod_{i=1}^n d_i -
\sum_{i=1}^n a_i c_i
\prod_{\stackrel{1\le j\le n,}{j\neq i}} d_j.
\end{equation}
Since $\rang A'= n+1$, then we have
$\rang A\ge n-1$.
If $\rang A =n-1$, then
$d_n=0$, implying
\[\det A'=-a_n c_n\prod_{i=1}^{n-1} d_i ;\]
otherwise
\[
\det A'=
\left( b d_n -
\sum_{i=1}^n a_i c_i d_n/d_i
\right)
\prod_{i=1}^{n-1}d_i .
\]
In both cases
$\prod_{i=1}^{n-1}d_i$ divides $\det A$. \qed
\end{enumerate}
 \end{pf}

Suppose
$A\in \mathcal{A}_n$,
$A'\in\bord(A)$.
From Table~\ref{SNFincid},
we see the following interesting facts:

\begin{itemize}
\item
The first $\circ$ in some $M_n$ corresponds to
$s=(1,0)$, $s'=(1,1,2)$. It is equivalent to following statement:
if
$A\in \mathcal{A}_{2,1}$ then $|\det A'|<2$.
\item if $A\in \mathcal{A}_{3,2}$, then $|\det A'|<3$.
\item if $A\in \mathcal{A}_4$, $\SNF(A)=(1,1,1,0)$, then $|\det A'|<5$.
\item if $A\in \mathcal{A}_4$, $\SNF(A)=(1,1,2,0)$, then $|\det A'|<4$.
\item if $A\in \mathcal{A}_5$, $\SNF(A)=(1,1,1,1,0)$,
then $|\det A'|\neq 7$.
\item if $A\in \mathcal{A}_5$, $\SNF(A)=(1,1,1,2,0)$,
then $|\det A'|\neq 6$.
\item if $A\in \mathcal{A}_5$, $\SNF(A)=(1,1,1,3,0)$,
then $|\det A'|\notin\{ 6,9\}$.
\item if $A\in \mathcal{A}_5$, $\SNF(A)=(1,1,1,2,2)$,
then $\SNF(A')\neq(1,1,1,1,4,0)$.
\item if $\SNF(A)=(1^{n-1},d_n)$ and  $\SNF(A')=(1^{n-1},d'_n,0)$
then $d'_n$ divides $d_n$ for all $n\le 7$.
\item if $\SNF(A)=s=(1^{n-1},d_n)\in\mathcal{S}_n$
and  $\SNF(A')=s'=(1^{n},d'_{n+1})\in\mathcal{S}_{n+1}$ then
 \begin{itemize}
 \item
 if $n\le 6$, then $m_{s,s'}=1$.
 \item
 if $n=7$, then $m_{s,s'}=1$ if and only if
\begin{eqnarray*}
(d_n,d_{n+1}')&\notin& \left\{(17,34), (7,39), (13,39), (1,42),\right.\\
 & & \left. (4,42), (6,42), (7,42), (13,42), (14,42)\right\}.
\end{eqnarray*}
 \item
 if $n=8$, then there are more exceptions to $m_{s,s'}=1$,
 but there is one exotic group of them:
 if $d_n=19$ then $d_{n+1}'$ must be divisible by $19$;
 $19$ is the only integer satisfying such a condition.
 \end{itemize}
\end{itemize}

\section{Determinant and SNF sets of \nj matrices of order $9$}
\label{snf9}

Determination of 
$\{|\det(A')| \mid A'\in\bord(A)\}$
is a simple operation, see the explanation following~\eqref{detobr}.
It was effectively performed
for all
$199727714$ matrices in $\mathcal{A}_8^\phi$;
merging these sets
$\mathcal{D}_9$ is obtained, see Table~\ref{brojanje}.

The similar idea --- determine ADV's, and only if necessary,
determine SNF's of the results of extension ---
is used to obtain
$\mathcal{S}_9$.
Suppose we know in advance the number $f_d$ of different SNF's
in $\mathcal{D}_9$
corresponding
to a given ADV $d>0$.
During the extension of
matrices from \ann,
the SNF's of extended matrices with the ADV $d$ are determined only
if the number of SNF's with ADV $d$ is still less
than $f_d$.
If we know only upper bound on $f_d$,
then the heuristic does not work ---
we have to determine SNF's of all matrices
with the ADV $d$.
Therefore, it is useful to determine $f_d$ for at least some $d>0$.

Denote by $p_n(r)$ the number of partitions of
$r$ into at most $n$ positive integers.
In order to determine the upper bound for $f_d$,
suppose first that $d$ is a prime power, $d=p^r$.
If $A\in\mathcal{A}_n$ and $|\det A|=d$, then $\SNF(A)$
is of the form
\[
(p^{x_1},p^{x_2},\ldots,p^{x_n}),\quad
0\le x_1\le x_2\le \cdots\le x_n,\quad
\sum_{i=1}^n x_i = r.
\]
The number of different exponent vectors $(x_1,x_2,\ldots,x_n)$ is equal to
$p_n(r)$. The values $p_n(m)$ are computed using the recurrence
(see for example~\cite{Balakrishnan}) $p_n(0)=1,$ $n\ge0$,
$p_0(r)=0$ for $r\ge 1$,
 and $p_n(r)=p_{n-1}(r)+p_n(r-n)$,
see Table~\ref{particije}.

\begin{table}
 \caption{The number of partitions of $r$ into at most $n$ positive
integers.}
 \label{particije}
{\small
\[
\begin{array}{||r@{\,}||@{\,}rrrrrrrr@{\,}||}
\hline\hline
~~~ r & 0 & 1 & 2 & 3 & 4 & 5 & 6 & 7 \\
n~~~ &   &   &   &   &   &   &   &   \\
\hline\hline
   0  & 1 & 0 & 0 & 0 & 0 & 0 & 0 & 0  \\
   1  & 1 & 1 & 1 & 1 & 1 & 1 & 1 & 1  \\
   2  & 1 & 1 & 2 & 2 & 3 & 3 & 4 & 4  \\
   3  & 1 & 1 & 2 & 3 & 4 & 5 & 7 & 8  \\
   4  & 1 & 1 & 2 & 3 & 5 & 6 & 9 &11  \\
   5  & 1 & 1 & 2 & 3 & 5 & 7 &10 &13  \\
   6  & 1 & 1 & 2 & 3 & 5 & 7 &11 &14  \\
   7  & 1 & 1 & 2 & 3 & 5 & 7 &11 &15  \\
\hline\hline
\end{array}
\]
}\end{table}

\begin{exmp}
If $d=8=2^3$ and $n=6$ we have $p_6(3)=3$;
$\SNF(A)$ is one of
$(1^5,8)$,
$(1^4,2,4)$ and
$(1^3,2,2,2)$. We see from Table~\ref{kat8} that all these SNF do
exist,
i.e. for each of them there exists some \nj matrix.
Another example
$d=3^2$, $n=6$, shows that
$p_6(2)=2$ is only an upper bound:
the SNF-class $(1^4,3,3)$ is empty.
\end{exmp}

More generally, if
$d=\prod_{i=1}^m p_i^{\alpha_i}$, where $p_i$ are different primes,
then
the upper bound on the number of different SNF's
with the ADV $d$ is
$\prod_{i=1}^m p_n(\alpha_i)$.

\begin{exmp}
If  $n=8$ and $d=36$,
then there are
$p_8(2)p_8(2)=4$ such  SNF's:
$(1^6,2,18)$,
$(1^7,36)$,
$(1^6,3,12)$,
$(1^6,6,6)$; all these SNF's are found in Table~\ref{kat8}.
\end{exmp}

In order to obtain a tighter upper bound for the
number of different SNF's, we have to include
somewhat more information.
If we further suppose that
$A'\in\bord(A)$  and
$\SNF(A)=s=(d_1,d_2,\ldots,d_{n})$,
then by Lemma~\ref{SNFogr} for some
$s'=(d_1',d_2',\ldots,d_{n+1}')$
the equality $\SNF(A')=s$ is
impossible.
For example, if $s$ contains $k$ ones,
then $s'$ contains at least $k$ ones
and if $\rang A\ge k+1$ then
$d'_{k+1}$ divides $d_{k+1}$.

Using these facts, the regular part of $\mathcal{S}_9$
was determined,
see Table~\ref{SNFdevet}.
The $\phi$\hyph representatives from
the chosen SNF-class of
$\mathcal{A}_8$
were extended computing determinants, and, if necessary,
determining SNF's.
The upper bounds for the number of different SNF's,
obtained by Lemma~\ref{SNFogr},
are rough for larger ADV values,
but the consequences are not dangerous,
because of the small number of extended matrices with the large ADV:
it is not hard to compute the SNF's of all of them.

To complete $\mathcal{S}_9$, it is necessary to determine the
singular part of $\mathcal{S}_9$.
If we would know that $F(8)$ is true,
then the set of singular SNF's of order $8$
would be equal to $\mathcal{S}_8$
(with each SNF extended by one zero, of course).
Not knowing a simple proof of $F(8)$, we proceed with a
shortened exhaustive proof.

The idea is to narrow the set of SNF-classes in $\mathcal{S}_8$,
the extension of which can lead to a  new
singular SNF of order $9$.
If
$\SNF(A)=(d_1,d_2,\ldots,d_8,0)$
for some $A\in\mathcal{A}_9$,
then
(because we know the set of SNF's of lower orders)
by Lemma~\ref{SNFogr} we can narrow the set SNF-classes,
containing $A$.
We obtain that
the only new possible SNF's are
the following SNF's of the rank $8$:
$(1^7,m,0),$ $m=44,45,48,56$ and $(1^6,2,28,0)$;
and the following SNF's
of the rank $7$:
$(1^6,20,0,0),$ $(1^6,24,0,0),$ $(1^6,32,0,0),$
$(1^5,2,12,0,0),$ $(1^5,2,16,0,0),$ $(1^5,4,8,0,0),$
$(1^4,2,2,8,0,0),$ $(1^4,2,4,4,0,0).$

The extension of which matrices gives the matrices with such SNF's?
For example, we know that the SNF $(1^7,44,0)$
can be obtained only by the extension of
a matrix in which
$44$ divides all minors of order $8$;
therefore $44$ also divides a nonsingular minor of order $8$;
hence the SNF of that minor could be only $(1^6,2,22)$.
Considering analogously the rest of listed SNF's of order $8$,
we obtain that matrices from $\mathcal{A}_{9,8}$,
with the SNF equal to some from the list above,
can be obtained only by the extension of
matrices from
$\mathcal{A}_8$ with the SNF $(1^6,2,22),$ $(1^6,2,24),$ $(1^6,3,15),$
$(1^5,2,2,12),$ or $(1^5,2,2,14).$

Analogously, we obtain that matrices from
$\mathcal{A}_{9,7}$ with one of the listed SNF's,
can be obtained only by double extension
of matrices from $\mathcal{A}_7$ with the SNF
$(1^5,2,10),$ $(1^4,2,2,6),$ or $(1^3,2,2,2,4).$
After the complete search through all
matrices that can be obtained by the extensions listed,
it is found that there are no new singular
SNF's of order $9$
i.e. that
$F(8)$ is also true.
That completes the determination of
$\mathcal{S}_9$.

In Table~\ref{SNFincidb} the part of the incidence matrix
$M_8$
is shown,
corresponding to regular matrices in
$\mathcal{S}_9$.
The table was obtained by extending
$\phi$\hyph representatives from
$\mathcal{A}_{8,7}$ and $\mathcal{A}_{8,8}$;
the singular extended matrices were ignored.

\section{The lower bounds for the first missing determinant,
$a_n$}\label{donjegr}

Denote by $f_n$ the $n$th Fibonacci number 
($f_1=f_2=1$ and $f_n=f_{n-1}+f_{n-2}$
for $n\ge 3$).
Paseman~\cite{paseman} shows that
$a_n\ge 2f_{n-1}$. 
We give the sketch of his proof,
and then we give
the sharper lower bounds for $a_n$,
$n\le 19$.

Consider the so called Fibonacci matrices
$F_n\in \mathcal{A}_n$
with the $(i,j)$ element equal to $1$
if and only if $j-i=-1,0,2,4,\ldots $;
$\det F_n=f_n$.
The cofactors corresponding to the first row of $F_n$ are
$f_{n-1}$, $f_{n-2}$, $-f_{n-3}$, $-f_{n-4}\ldots,-f_{1}$.
Consider the matrix $U\in\bord(F_n)$,

\[
U=
\left[\begin{array}{rr}
  F_n & y \\
  x & b
\end{array}\right],
\]
where $x=[x_1\ x_2\ \ldots\ x_{n-1}]$, 
$y=[y_1\ y_2\ \ldots\ y_{n-1}]^T$.
Let $y_1=1$, $y_2=y_3=\cdots=y_n=0$
and $x_1=x_2=0$. 
Then from ~\eqref{detobr} we have
\[
\det U=\sum_{i=1}^{n-2}x_{n+1-i}f_i+ bf_n.
\]
Therefore, each integer from $[0,2f_n-1]$ is determinant of some
$U\in\bord(F_n)$, and $a_n\ge 2f_{n-1}$.

In order to prove that $a_n\ge m$,
one can give a list of matrices from $\mathcal{A}_{n-1}$,
such that determinants of their extensions cover $[1,m-1]$.
The proof verification then includes the procedure of
finding determinants of all extensions of a given matrix.
Still, such a list is essentially more compact than the list of matrices from
$\mathcal{A}_{n}$,
with determinants covering $[1,m-1]$.

Denote by $a_A$ 
the minimal integer not in
$\cup\{|\det B| \mid B\in\bord(A)\}$, the
''extension spectrum'' of $A\in\mathcal{A}_{n}$.
In this context, 
the matrices $A$ with high $a_A$ are of special interest.
If $a_A>1$ and $\SNF(A)=(d_1,d_2,\ldots,d_n)$,
then
$d_1=d_2=\cdots=d_{n-1}=1$, because determinants of all extensions of
$A$ are divisible by $d_{n-1}$, see~\eqref{z}.

In order to find lower bounds for some $a_n$,
one can start from a well chosen set
$\mathcal{B}_{n-1}\subset
\mathcal{A}_{n-1}$,
and then to find ADV's of all extended matrices.
If $m$ is the smallest number not equal to some
of these ADV's, then
$a_n\ge m$.
Afterwards, some subset of extended matrices with different SNF's is taken
to be the set
$\mathcal{B}_{n}$,
and the next iteration can be started.

The starting set
$\mathcal{B}_{9}$
was constructed in the following way.
From each SNF-class in
$\mathcal{A}_{8}$
a number of matrices is taken,
with different numbers of $\pi$ representatives in
their $\phi$\hyph classes.
Extending these matrices, a set of matrices with different
SNF's is obtained, but without any matrix
with the SNF
$(1^8,97)$.
By adding one such matrix, the set
$\mathcal{B}_{9}$ is completed.
The sets
$\mathcal{B}_{10}$,
$\mathcal{B}_{11}$ and
$\mathcal{B}_{12}$ are generated iteratively,
as explained above.
At the end, the ADV's of all matrices
obtained by extending the matrices in
$\mathcal{B}_{12}$ are determined.
The resulting lower bounds are
$a_{10}\ge 259$,
$a_{11}\ge 739$,
$a_{12}\ge 2107$,
$a_{13}\ge 6157$.

For $n>13$ we used an alternative heuristic, described by Algorithm~\ref{heuristika}.

\begin{alg}\label{heuristika}
Heuristic to find lower bound for $a_{n+1}$.
\begin{tabbing}
 \quad \= \quad \= \quad \= \quad \= \quad \= \quad \=\quad  \kill
  \INPUT  $\mathcal{L}_n\subset \mathcal{A}_n$, list of matrices to be extended.  \\
  \OUTPUT  lower bound for $a_{n+1}$, and list $\mathcal{L}_{n+1}\subset \mathcal{A}_{n+1}$\\
  \>  of ''promising'' matrices for the following iteration.  \\
  \{ Initialization: \} \\
  $first0\la 1$; \{ the first integer not ''covered'' by ADV's \} \\
  $dmax\la 1$;  \{ the largest ADV found until now \} \\
  $\mathcal{L}_{n+1}\la \emptyset$; \{ output list \} \\
  \FOR all $A\in\mathcal{L}_n$ \DO \\
  \>\{ Consider the extensions $A'=\left[\begin{array}{rr} A & y \\  x & b\end{array}\right]$\} \\
  \> Compute $\det A$ and $B=[B_{ij}]=\mathrm{adj}\, A$; \{ transposed cofactor matrix of $A$\} \\
  \>  \FOR all $y\in\{0,1\}^n$ \DO \\
  \>\> \{ the next linear combination of rows of $B$ \} \\
  \>\> determine the coefficients of the linear combination \\
  \>\>\> $-b\det A + \sum_{i=1}^n x_i \left( \sum_{j=1}^n y_jB_{ij}\right)$ \\
  \>\>\> and the sums $s^+$, $s^-$ of its positive and negative members; \\
  \>\> \IF $\max\{-s^-,s^+\}\ge first0$ \THEN \{''poor'' linear combinations are skipped\} \\
  \>\>\> \FOR all $(x,b)\in\{0,1\}^{n+1}$ \DO \\
  \>\>\>\> compute $\det A'$; \{ by one addition only, using Gray code \} \\
  \>\>\>\> \IF $|\det A'|=first0$ \THEN \\
  \>\>\>\>\> update $first0$; \\
  \>\>\>\> \IF $|\det A'|>dmax$ \THEN \\
  \>\>\>\>\> $dmax\la |\det A'|$; \\
  \>\>\>\> \IF $|\det A'|>0.9\, dmax$ \THEN \\
  \>\>\>\>\> append $A'$ to $\mathcal{L}_{n+1}$;\\
\end{tabbing}
\end{alg}

Elimination of ''poor'' linear combinations is a powerful heuristics
if the matrices with the high extension spectra are placed in the beginning of
$\mathcal{L}_{n}$.
The major part of linear combinations is skipped after only a few first matrices in
$\mathcal{L}_{n}$,
reducing the extension complexity roughly to $O(n2^n)$ (instead of $O(4^n)$).
In Table~\ref{prve0} for $10\le n\le 19$ we give
\begin{itemize}
  \item lower bound for $a_n$,
  \item $|\mathcal{L}_{n-1}|$, the number of extended matrices,
  \item a matrix $A_{n-1}$ with the highest extension spectrum found in $|\mathcal{A}_{n-1}|$,
  \item extension spectrum and determinant of $A_{n-1}$.
\end{itemize}
Complete lists of matrices, whose extension determinants
prove these lower bounds,
can be fount at
http://www.matf.bg.ac.yu/~ezivkovm/01matrices.htm.

\section{Counting \nj matrices with the maximum determinant}
\label{brojevi}

Using the classification of \ann, it is not hard to compute
the number $c_n$~\cite[Sequences A051752]{Slo}
 of matrices in \an with the maximal determinant
$d_n$ (i.e. $1/2$ of the number of matrices with the ADV $d_n$)
for $n\le 9$.

The first $8$ members of the sequence $c_n$
are found in Table~\ref{kat8};
the number $c_8=195955200$ is new.

In order to determine $c_9$,
from Table~\ref{SNFincid} we see that the matrix from
$\mathcal{A}_9$
with the ADV $144$ could be obtained only by 
extending matrices from
$\mathcal{A}_8$
with the SNF
$(1^5,2,2,6)$ or
$(1^5,2,2,12)$.
After the extension of these two SNF-classes,
it turned out that there is a unique $\phi$\hyph class with the ADV
$144$ --- the class with the representative
$(\text{F,33,C3,FC,155,15A,166,196,1A9})$.
Half of the number of matrices in that $\phi$\hyph class is
$c_9=13716864000$.
It is interesting that for all $n\le 9$
there is a unique
$\phi$\hyph class
with the maximal ADV.

\section{Acknowledgement}

I am greatly indebted to the anonymous referee whose comments helped
to improve the exposition.

\bibliographystyle{amsplain}

\appendix

\section{Large tables}

\begin{table}
 \caption{$\pi$\hyph representatives of \nj matrices of order $3$.}
\label{pi3}
{\small

\begin{tabular}{||rrr|r||}
                                         \hline\hline
                                         \multicolumn{3}{||c|}{SNF} & size \\
  0&   0&   0&    1\\
                                         \hline\hline
  0&   0&   0&    1\\
                                         \hline\hline
                                         \multicolumn{4}{c}{} \\
                                         \multicolumn{4}{c}{} \\
                                         \multicolumn{4}{c}{} \\
                                         \multicolumn{4}{c}{} \\
                                         \multicolumn{4}{c}{} \\
                                         \hline\hline
                                         \multicolumn{3}{||c|}{SNF} & size \\
  1&   0&   0&   49\\
                                         \hline\hline
  0&   0&   1&    9\\
  0&   0&   7&    3\\
  1&   1&   1&    3\\
  7&   7&   7&    1\\
                                         \hline
  0&   0&   3&    9\\
  3&   3&   3&    3\\
                                         \hline
  0&   1&   1&    9\\
  0&   7&   7&    3\\
                                         \hline
  0&   3&   3&    9\\
                                         \hline\hline
                                         \end{tabular}\qquad
                                         \begin{tabular}{||rrr|r||}
                                         \hline\hline
                                         \multicolumn{3}{||c|}{SNF} & size \\
  1&   1&   0&  288\\
                                         \hline\hline
  0&   1&   2&   18\\
  0&   1&   7&   18\\
  1&   1&   3&   18\\
  1&   1&   6&    9\\
  3&   7&   7&    9\\
                                         \hline
  0&   1&   3&   36\\
  0&   1&   6&   18\\
  0&   3&   7&   18\\
  1&   1&   2&   18\\
  1&   3&   3&   18\\
  1&   1&   7&    9\\
  3&   3&   7&    9\\
  1&   6&   6&    9\\
  1&   7&   7&    9\\
                                         \hline
  0&   3&   5&   18\\
  3&   3&   5&   18\\
                                         \hline
  1&   2&   3&   18\\
  1&   6&   7&   18\\
                                         \hline\hline
                                         \end{tabular}\qquad
                                         \begin{tabular}{||rrr|r||}
                                         \hline\hline
                                         \multicolumn{3}{||c|}{SNF} & size \\
  1&   1&   1&  168\\
                                         \hline\hline
  1&   2&   4&    6\\
  1&   2&   7&   18\\
  1&   3&   5&   18\\
  1&   3&   6&   36\\
  3&   5&   7&   18\\
                                         \hline
  1&   2&   5&   36\\
  1&   3&   7&   36\\
                                         \hline\hline
                                         \multicolumn{4}{c}{} \\
                                         \multicolumn{4}{c}{} \\
                                         \multicolumn{4}{c}{} \\
                                         \multicolumn{4}{c}{} \\
                                         \multicolumn{4}{c}{} \\
                                         \multicolumn{4}{c}{} \\
                                         \multicolumn{4}{c}{} \\
                                         \hline\hline
                                         \multicolumn{3}{||c|}{SNF} & size \\
  1&   1&   2&    6\\
                                         \hline\hline
  3&   5&   6&    6\\
                                         \hline\hline
                                         \end{tabular}

}
\end{table}

\begin{table}
 \caption{$\phi$\hyph representatives of \nj matrices of order $4$.}
\label{phi4}
{\small
\begin{tabular}{||rrrr|r||}
\hline\hline
\multicolumn{4}{||c|}{SNF} & size                                       \\
  0& 0&  0&  0 & 1                                     \\
\hline\hline
  0& 0&  0&  0 &   1                                     \\
\hline\hline
\multicolumn{4}{||c|}{SNF} & size                                      \\
  1& 0&  0&  0 & 225                                     \\
\hline
  0& 0&  0&  1 &  25                                     \\
  0& 0&  0&  3 &  50                                     \\
  0& 0&  1&  1 &  50                                     \\
  0& 0&  3&  3 & 100                                     \\
\hline\hline
\multicolumn{4}{||c|}{SNF} & size                                      \\
  1& 1&  0&  0 &6750                                     \\
\hline\hline
  0& 0&  1&  2 & 200                                     \\
  0& 0&  1&  3 & 400                                     \\
  0& 0&  1&  6 & 600                                     \\
  0& 0&  3&  5 & 600                                     \\
  0& 0&  3&  7 & 300                                     \\
  0& 1&  1&  2 & 600                                     \\
  0& 1&  1&  6 & 450                                     \\
  0& 1&  2&  3 & 600                                     \\
  0& 1&  3&  3 & 300                                     \\
  0& 1&  6&  6 & 900                                     \\
  0& 1&  6&  7 & 900                                     \\
  0& 3&  3&  5 & 900                                     \\
\hline\hline
\end{tabular}\quad
\raisebox{10ex}{
\begin{tabular}{||rrrr|r||}
\hline\hline
\multicolumn{4}{||c|}{SNF} & size                                      \\
  1& 1&  1&  0 &35400                                     \\
\hline\hline
  0& 1&  2&  4 & 600                                     \\
  0& 1&  2&  5 &3600                                     \\
  0& 1&  2&  7 &1800                                     \\
  0& 1&  3&  5 &1800                                     \\
  0& 1&  3&  6 &3600                                     \\
  0& 1&  3&  7 &3600                                     \\
  0& 1&  6&  A &3600                                     \\
  0& 3&  5&  7 &1800                                     \\
  0& 3&  5&  9 &1200                                     \\
  0& 3&  5&  A &3600                                     \\
  1& 2&  3&  4 &3600                                     \\
  1& 2&  4&  7 &1200                                     \\
  1& 2&  5&  6 &3600                                     \\
  1& 6&  7&  A &1800                                     \\
\hline\hline
\end{tabular}
}\quad
\begin{tabular}{||rrrr|r||}
\hline\hline
\multicolumn{4}{||c|}{SNF} & size                                      \\
  1& 1&  2&  0 & 600                                     \\
\hline\hline
 3&  5&  6 &  0& 600                                     \\
\hline\hline
\multicolumn{5}{c}{}  \\
\multicolumn{5}{c}{}  \\
\hline\hline
\multicolumn{4}{||c|}{SNF} & size                                      \\
  1& 1&  1&  1 &20040                                     \\
\hline\hline
  1& 2&  4&  8 & 600                                     \\
  1& 2&  4&  9 &7200                                     \\
  1& 2&  5&  A &1440                                     \\
  1& 2&  5&  B &7200                                     \\
  1& 2&  7&  B &3600                                     \\
\hline\hline
\multicolumn{5}{c}{}  \\
\hline\hline
\multicolumn{4}{||c|}{SNF} & size                                      \\
  1& 1&  1&  2 &2400                                     \\
\hline\hline
  1& 6&  A&  C &2400                                     \\
\hline\hline
\multicolumn{5}{c}{}  \\
\multicolumn{5}{c}{}  \\
\hline\hline
\multicolumn{4}{||c|}{SNF} & size                                      \\
  1& 1&  1&  3 &120                                     \\
\hline\hline
  3& 5&  9&  E & 120                                     \\
\hline\hline
\end{tabular}
}
\end{table}

\begin{table}
 \caption{The representatives and the sizes of SNF-classes in \ann,
$n\le 8$.} \label{kat8}
{\small
\renewcommand{\arraystretch}{1.0}
\[
\begin{array}{||r|r|r||r|r|r||r||}
\hline\hline
\multicolumn{3}{||c||}{\text{$\mathcal{A}_1$}}&
\multicolumn{3}{c||}{\text{The number of}}&
\multicolumn{1}{c||}{\text{The SNF-class}} \\
\cline{1-3}
\text{~}&\det& \multicolumn{1}{c||}{\text{SNF}}&\text{matrices } &\text{$\pi$-classes}&\text{$\phi$-classes}& \multicolumn{1}{c||}{\text{representative}}     \\
\hline\hline
    0&      0& 0&  1&  1&  1&  0\\
\hline
    1&      1& 1&  1&  1&  1&  1\\
\hline\hline
\multicolumn{3}{||c||}{\text{Total:}}&2&  2&  2&   \\
\hline\hline
\end{array}
\]
}

{\small
\renewcommand{\arraystretch}{1.0}
\[
\begin{array}{||r|r|r@{\rz}r||r|r|r||r@{\rz}r||}
\hline\hline
\multicolumn{4}{||c||}{\text{$\mathcal{A}_2$}}&
\multicolumn{3}{c||}{\text{The number of}}&
\multicolumn{2}{c||}{\text{The SNF-class}} \\
\cline{1-4}
\text{~}&\det& \multicolumn{2}{c||}{\text{SNF}}&\text{matrices } &\text{$\pi$-classes}&\text{$\phi$-classes}& \multicolumn{2}{c||}{\text{representative}}     \\
\hline\hline
    0&      0&0& 0&  1&  1&  1&  0&  0\\
\hline
    1&      0&0& 1&  9&  4&  1&  0&  1\\
\hline
    2&      1&1& 1&  6&  2&  1&  1&  2\\
\hline\hline
\multicolumn{4}{||c||}{\text{Total:}}&16& 7&  3&   &   \\
\hline\hline
\end{array}
\]
}

{\small
\renewcommand{\arraystretch}{1.0}
\[
\begin{array}{||r|r|r@{\rz}r@{\rz}r||r|r|r||r@{\rz}r@{\rz}r||}
\hline\hline
\multicolumn{5}{||c||}{\text{$\mathcal{A}_3$}}&
\multicolumn{3}{c||}{\text{The number of}}&
\multicolumn{3}{c||}{\text{The SNF-class}} \\
\cline{1-5}
\text{~}&\det& \multicolumn{3}{c||}{\text{SNF}}&\text{matrices } &\text{$\pi$-classes}&\text{$\phi$-classes}& \multicolumn{3}{c||}{\text{representative}}     \\
\hline\hline
    0&      0&0&0& 0&    1&  1&  1&  0&  0&  0\\
\hline
    1&      0&0&0& 1&   49&  9&  4&  0&  0&  1\\
\hline
    2&      0&0&1& 1&  288& 18&  4&  0&  1&  2\\
\hline
    3&      1&1&1& 1&  168&  7&  2&  1&  2&  4\\
    4&      2&1&1& 2&    6&  1&  1&  3&  5&  6\\
\hline    \hline
\multicolumn{5}{||c||}{\text{Total:}}&512& 36& 12&   &   &   \\
\hline\hline
\end{array}
\]
}

{\small
\renewcommand{\arraystretch}{1.0}
\[
\begin{array}{||r|r|r@{\rz}r@{\rz}r@{\rz}r||r|r|r||r@{\rz}r@{\rz}r@{\rz}r||}
\hline\hline
\multicolumn{6}{||c||}{\text{$\mathcal{A}_4$}}&
\multicolumn{3}{c||}{\text{The number of}}&
\multicolumn{4}{c||}{\text{The SNF-class}} \\
\cline{1-6}
\text{~}&\det& \multicolumn{4}{c||}{\text{SNF}}&\text{matrices } &\text{$\pi$-classes}&\text{$\phi$-classes}& \multicolumn{4}{c||}{\text{representative}}     \\
\hline\hline
    0&      0&0&0&0& 0&      1&   1&  1&\text{ 0}&\text{ 0}&\text{ 0}&\text{ 0}\\
\hline
    1&      0&0&0&0& 1&    225&  16&  4&\text{ 0}&\text{ 0}&\text{ 0}&\text{ 1}\\
\hline
    2&      0&0&0&1& 1&   6750&  84& 12&\text{ 0}&\text{ 0}&\text{ 1}&\text{ 2}\\
\hline
    3&      0&0&1&1& 1&  35400& 150& 14&\text{ 0}&\text{ 1}&\text{ 2}&\text{ 4}\\
    4&      0&0&1&1& 2&    600&   5&  1&\text{ 0}&\text{ 3}&\text{ 5}&\text{ 6}\\
\hline
    5&      1&1&1&1& 1&  20040&  49&  5&\text{ 1}&\text{ 2}&\text{ 4}&\text{ 8}\\
    6&      2&1&1&1& 2&   2400&  10&  1&\text{ 1}&\text{ 6}&\text{ A}&\text{ C}\\
    7&      3&1&1&1& 3&    120&   2&  1&\text{ 3}&\text{ 5}&\text{ 9}&\text{ E}\\
\hline\hline
\multicolumn{6}{||c||}{\text{Total:}}&65536& 317& 39&   &   &   &   \\
\hline\hline
\end{array}
\]
}
\end{table}

\renewcommand{\arraystretch}{1.0}
\begin{table}
\addtocounter{table}{-1}
 \caption{Continued}
      \centering
{\small
\[
\begin{array}{||r|r|r@{\rz}r@{\rz}r@{\rz}r@{\rz}r||r|r|r||r@{\rz}r@{\rz}r@{\rz}r@{\rz}r||}
\hline\hline
\multicolumn{7}{||c||}{\text{$\mathcal{A}_5$}}&
\multicolumn{3}{c||}{\text{The number of}}&
\multicolumn{5}{c||}{\text{The SNF-class}} \\
\cline{1-7}
\text{~}&\det& \multicolumn{5}{c||}{\text{SNF}}&\text{matrices } &\text{$\pi$-classes}&\text{$\phi$-classes}& \multicolumn{5}{c||}{\text{representative}}     \\
\hline\hline
    0&      0&0&0&0&0& 0&         1&    1&   1&\text{ 0}&\text{ 0}&\text{ 0}&\text{ 0}&\text{ 0}\\
\hline
    1&      0&0&0&0&0& 1&       961&   25&   9&\text{ 0}&\text{ 0}&\text{ 0}&\text{ 0}&\text{ 1}\\
\hline
    2&      0&0&0&0&1& 1&    118800&  260&  37&\text{ 0}&\text{ 0}&\text{ 0}&\text{ 1}&\text{ 2}\\
\hline
    3&      0&0&0&1&1& 1&   3134400& 1346& 113&\text{ 0}&\text{ 0}&\text{ 1}&\text{ 2}&\text{ 4}\\
    4&      0&0&0&1&1& 2&     25350&   25&   5&\text{ 0}&\text{ 0}&\text{ 3}&\text{ 5}&\text{ 6}\\
\hline
    5&      0&0&1&1&1& 1&  16853400& 2589& 141&\text{ 0}&\text{ 1}&\text{ 2}&\text{ 4}&\text{ 8}\\
    6&      0&0&1&1&1& 2&    880200&  210&  17&\text{ 0}&\text{ 1}&\text{ 6}&\text{ A}&\text{ C}\\
    7&      0&0&1&1&1& 3&     27000&   15&   2&\text{ 0}&\text{ 3}&\text{ 5}&\text{ 9}&\text{ E}\\
\hline
    8&      1&1&1&1&1& 1&   9702720&  831&  39&\text{ 1}&\text{ 2}&\text{ 4}&\text{ 8}&\text{10}\\
    9&      2&1&1&1&1& 2&   2427840&  254&  15&\text{ 1}&\text{ 2}&\text{ C}&\text{14}&\text{18}\\
   10&      3&1&1&1&1& 3&    289440&   51&   5&\text{ 1}&\text{ 6}&\text{ A}&\text{12}&\text{1C}\\
   11&      4&1&1&1&1& 4&     65520&   12&   2&\text{ 3}&\text{ 5}&\text{ 9}&\text{11}&\text{1E}\\
   12&      5&1&1&1&1& 5&      7200&    3&   1&\text{ 3}&\text{ 5}&\text{ E}&\text{16}&\text{19}\\
   13&      4&1&1&1&2& 2&     21600&    2&   1&\text{ 3}&\text{ C}&\text{15}&\text{16}&\text{19}\\
\hline\hline
\multicolumn{7}{||c||}{\text{Total:}}&33554432& 5624& 388&   &   &   &   &   \\
\hline\hline
\end{array}
\]
%
\renewcommand{\arraystretch}{0.8}
\[
\begin{array}{||r|r|r@{\rz}r@{\rz}r@{\rz}r@{\rz}r@{\rz}r||r|r|r||r@{\rz}r@{\rz}r@{\rz}r@{\rz}r@{\rz}r||}
\hline\hline
\multicolumn{8}{||c||}{\text{$\mathcal{A}_6$}}&
\multicolumn{3}{c||}{\text{The number of}}&
\multicolumn{6}{c||}{\text{The SNF-class}} \\
\cline{1-8}
\text{~}&\det& \multicolumn{6}{c||}{\text{SNF}}&\text{matrices } &\text{$\pi$-classes}&\text{$\phi$-classes}& \multicolumn{6}{c||}{\text{representative}}     \\
\hline\hline
    0&    0&0&0&0&0&0& 0&            1&      1&    1&\text{ 0}&\text{ 0}&\text{ 0}&\text{ 0}&\text{ 0}&\text{ 0}\\
\hline
    1&    0&0&0&0&0&0& 1&         3969&     36&    9&\text{ 0}&\text{ 0}&\text{ 0}&\text{ 0}&\text{ 0}&\text{ 1}\\
\hline
    2&    0&0&0&0&0&1& 1&      1807806&    660&   76&\text{ 0}&\text{ 0}&\text{ 0}&\text{ 0}&\text{ 1}&\text{ 2}\\
\hline
    3&    0&0&0&0&1&1& 1&    189336000&   7586&  472&\text{ 0}&\text{ 0}&\text{ 0}&\text{ 1}&\text{ 2}&\text{ 4}\\
    4&    0&0&0&0&1&1& 2&       735000&     86&   10&\text{ 0}&\text{ 0}&\text{ 0}&\text{ 3}&\text{ 5}&\text{ 6}\\
\hline
    5&    0&0&0&1&1&1& 1&   5168108400&  47605& 1913&\text{ 0}&\text{ 0}&\text{ 1}&\text{ 2}&\text{ 4}&\text{ 8}\\
\hline\hline
\end{array}
\]
}
\end{table}

\renewcommand{\arraystretch}{1.0}

\begin{table}
\addtocounter{table}{-1}
 \caption{Continued}
      \centering
{\small
\[
\begin{array}{||r|r|r@{\rz}r@{\rz}r@{\rz}r@{\rz}r@{\rz}r@{\rz}r||r|r|r||r@{\rz}r@{\rz}r@{\rz}r@{\rz}r@{\rz}r@{\rz}r||}
\hline\hline
\multicolumn{9}{||c||}{\text{$\mathcal{A}_7$}}&
\multicolumn{3}{c||}{\text{The number of}}&
\multicolumn{7}{c||}{\text{The SNF-class}} \\
\cline{1-9}
\text{~}&\det& \multicolumn{7}{c||}{\text{SNF}}&\text{matrices } &\text{$\pi$-classes}&\text{$\phi$-classes}& \multicolumn{7}{c||}{\text{representative}}     \\
\hline\hline
    0&    0&0&0&0&0&0&0& 0&              1&        1&      1&\text{ 0}&\text{ 0}&\text{ 0}&\text{ 0}&\text{ 0}&\text{ 0}&\text{ 0}\\
\hline
    1&    0&0&0&0&0&0&0& 1&          16129&       49&     16&\text{ 0}&\text{ 0}&\text{ 0}&\text{ 0}&\text{ 0}&\text{ 0}&\text{ 1}\\
\hline
    2&    0&0&0&0&0&0&1& 1&       25316928&     1428&    170&\text{ 0}&\text{ 0}&\text{ 0}&\text{ 0}&\text{ 0}&\text{ 1}&\text{ 2}\\
\hline
    3&    0&0&0&0&0&1&1& 1&     9254300328&    31994&   1908&\text{ 0}&\text{ 0}&\text{ 0}&\text{ 0}&\text{ 1}&\text{ 2}&\text{ 4}\\
    4&    0&0&0&0&0&1&1& 2&       17360406&      246&     34&\text{ 0}&\text{ 0}&\text{ 0}&\text{ 0}&\text{ 3}&\text{ 5}&\text{ 6}\\
\hline
    5&    0&0&0&0&1&1&1& 1&   989588124000&   501563&  17596&\text{ 0}&\text{ 0}&\text{ 0}&\text{ 1}&\text{ 2}&\text{ 4}&\text{ 8}\\
\hline\hline
\end{array}
\]
}
\end{table}

\begin{table}
\addtocounter{table}{-1}
 \caption{Continued}
      \centering
{\small
\[
\begin{array}{||r|r|r@{\rz}r@{\rz}r@{\rz}r@{\rz}r@{\rz}r@{\rz}r||r|r|r||r@{\rz}r@{\rz}r@{\rz}r@{\rz}r@{\rz}r@{\rz}r||}
\hline\hline
\multicolumn{9}{||c||}{\text{$\mathcal{A}_7$}}&
\multicolumn{3}{c||}{\text{The number of}}&
\multicolumn{7}{c||}{\text{The SNF-class}} \\
\cline{1-9}
\text{~}&\det& \multicolumn{7}{c||}{\text{SNF}}&\text{matrices } &\text{$\pi$-classes}&\text{$\phi$-classes}& \multicolumn{7}{c||}{\text{representative}}     \\
\hline\hline
   26&    1&1&1&1&1&1&1& 1&135491563468800&  5593528&  91764&\text{ 1}&\text{ 2}&\text{ 4}&\text{ 8}&\text{10}&\text{20}&\text{40}\\
   27&    2&1&1&1&1&1&1& 2& 83220427382400&  3493129&  58179&\text{ 1}&\text{ 2}&\text{ 4}&\text{ 8}&\text{30}&\text{50}&\text{60}\\
   28&    3&1&1&1&1&1&1& 3& 23436399974400&  1020752&  17707&\text{ 1}&\text{ 2}&\text{ 4}&\text{18}&\text{28}&\text{48}&\text{70}\\
   29&    4&1&1&1&1&1&1& 4& 13285672243200&   581948&  10189&\text{ 1}&\text{ 2}&\text{ C}&\text{14}&\text{24}&\text{44}&\text{78}\\
   30&    5&1&1&1&1&1&1& 5&  3754520017920&   172714&   3169&\text{ 1}&\text{ 2}&\text{ C}&\text{14}&\text{38}&\text{58}&\text{64}\\
   31&    6&1&1&1&1&1&1& 6&  4201407745920&   185688&   3320&\text{ 1}&\text{ 6}&\text{ A}&\text{12}&\text{2C}&\text{5C}&\text{62}\\
\hline\hline
\end{array}
\]
}
\end{table}


\begin{table}
\addtocounter{table}{-1}
 \caption{Continued}
      \centering
{\small
\[
\begin{array}{||r|r|r@{\rz}r@{\rz}r@{\rz}r@{\rz}r@{\rz}r@{\rz}r@{\rz}r||r|r|r||r@{\rz}r@{\rz}r@{\rz}r@{\rz}r@{\rz}r@{\rz}r@{\rz}r||}
\hline\hline
\multicolumn{10}{||c||}{\text{$\mathcal{A}_8$}}&
\multicolumn{3}{c||}{\text{The number of}}&
\multicolumn{8}{c||}{\text{The SNF-class}} \\
\cline{1-10}
\text{~}&\det& \multicolumn{8}{c||}{\text{SNF}}&\text{matrices } &\text{$\pi$-classes}&\text{$\phi$-classes}& \multicolumn{8}{c||}{\text{representative}}     \\
\hline\hline
    0&    0&0&0&0&0&0&0&0& 0&                    1&           1&         1&\text{  0}&\text{ 0}&\text{ 0}&\text{ 0}&\text{ 0}&\text{ 0}&\text{ 0}&\text{ 0}\\
\hline
    1&    0&0&0&0&0&0&0&0& 1&                65025&          64&        16&\text{  0}&\text{ 0}&\text{ 0}&\text{ 0}&\text{ 0}&\text{ 0}&\text{ 0}&\text{ 1}\\
\hline
    2&    0&0&0&0&0&0&0&1& 1&            336954750&        2800&       295&\text{  0}&\text{ 0}&\text{ 0}&\text{ 0}&\text{ 0}&\text{ 0}&\text{ 1}&\text{ 2}\\
\hline
    3&    0&0&0&0&0&0&1&1& 1&         396683821800&      110064&      5758&\text{  0}&\text{ 0}&\text{ 0}&\text{ 0}&\text{ 0}&\text{ 1}&\text{ 2}&\text{ 4}\\
    4&    0&0&0&0&0&0&1&1& 2&            362237400&         596&        52&\text{  0}&\text{ 0}&\text{ 0}&\text{ 0}&\text{ 0}&\text{ 3}&\text{ 5}&\text{ 6}\\
\hline
    5&    0&0&0&0&0&1&1&1& 1&      144191294561160&     3696215&    114651&\text{  0}&\text{ 0}&\text{ 0}&\text{ 0}&\text{ 1}&\text{ 2}&\text{ 4}&\text{ 8}\\
\hline\hline
\end{array}
\]
}
\end{table}

\begin{table}
\addtocounter{table}{-1}
 \caption{Continued}
      \centering
{\small
\[
\begin{array}{||r|r|r@{\rz}r@{\rz}r@{\rz}r@{\rz}r@{\rz}r@{\rz}r@{\rz}r||r|r|r||r@{\rz}r@{\rz}r@{\rz}r@{\rz}r@{\rz}r@{\rz}r@{\rz}r||}
\hline\hline
\multicolumn{10}{||c||}{\text{$\mathcal{A}_8$}}&
\multicolumn{3}{c||}{\text{The number of}}&
\multicolumn{8}{c||}{\text{The SNF-class}} \\
\cline{1-10}
\text{~}&\det& \multicolumn{8}{c||}{\text{SNF}}&\text{matrices } &\text{$\pi$-classes}&\text{$\phi$-classes}& \multicolumn{8}{c||}{\text{representative}}     \\
\hline\hline
   44&    0&0&1&1&1&1&1&2& 2&    25998918420787200&    31764328&    476566&\text{  0}&\text{ 1}&\text{ 2}&\text{ C}&\text{30}&\text{54}&\text{58}&\text{64}\\
   45&    0&0&1&1&1&1&1&2& 4&      904945412044800&     1434616&     23722&\text{  0}&\text{ 1}&\text{ 6}&\text{18}&\text{2A}&\text{34}&\text{4C}&\text{52}\\
   46&    0&0&1&1&1&1&1&2& 6&       38868105830400&       70882&      1251&\text{  0}&\text{ 3}&\text{ 5}&\text{19}&\text{29}&\text{36}&\text{4E}&\text{51}\\
   47&    0&0&1&1&1&1&1&2& 8&        1777705574400&        4348&       114&\text{  0}&\text{ 3}&\text{ C}&\text{31}&\text{55}&\text{5A}&\text{66}&\text{69}\\
   48&    0&0&1&1&1&1&1&2&10&         175575859200&         408&        10&\text{  0}&\text{ 7}&\text{19}&\text{2A}&\text{34}&\text{4C}&\text{52}&\text{63}\\
   49&    0&0&1&1&1&1&1&3& 3&       70619902617600&      117219&      1972&\text{  0}&\text{ 3}&\text{ 5}&\text{18}&\text{28}&\text{49}&\text{4E}&\text{71}\\
\hline\hline
\end{array}
\]
}
\end{table}

\begin{table}
\addtocounter{table}{-1}
 \caption{Continued}
      \centering
{\small
\[
\begin{array}{||r@{\rz}r|r@{\rz}r@{\rz}r@{\rz}r@{\rz}r@{\rz}r@{\rz}r@{\rz}r||r|r|r||r@{\rz}r@{\rz}r@{\rz}r@{\rz}r@{\rz}r@{\rz}r@{\rz}r||}
\hline\hline
\multicolumn{10}{||c||}{\text{$\mathcal{A}_8$}}&
\multicolumn{3}{c||}{\text{The number of}}&
\multicolumn{8}{c||}{\text{The SNF-class}} \\
\cline{1-10}
\text{~}&\det& \multicolumn{8}{c||}{\text{SNF}}&\text{matrices } &\text{$\pi$-classes}&\text{$\phi$-classes}& \multicolumn{8}{c||}{\text{representative}}     \\
\hline\hline
   87&  32&1&1&1&1&1&1&1&32&        5423648025600&        3953&        68&\text{  3}&\text{ C}&\text{31}&\text{56}&\text{6A}&\text{9A}&\text{A6}&\text{C5}\\
   88&  33&1&1&1&1&1&1&1&33&        2806470374400&        1866&        27&\text{  3}&\text{ D}&\text{1E}&\text{35}&\text{66}&\text{79}&\text{AA}&\text{D5}\\
   89&  34&1&1&1&1&1&1&1&34&        2831160729600&        1815&        23&\text{  3}&\text{ D}&\text{31}&\text{54}&\text{6A}&\text{9A}&\text{A6}&\text{C9}\\
   90&  35&1&1&1&1&1&1&1&35&         757170892800&         491&         7&\text{  3}&\text{1D}&\text{2E}&\text{56}&\text{69}&\text{9A}&\text{B1}&\text{CD}\\
   91&  36&1&1&1&1&1&1&1&36&        1327792435200&         873&        13&\text{  3}&\text{ D}&\text{35}&\text{56}&\text{69}&\text{99}&\text{AE}&\text{C5}\\
   92&  37&1&1&1&1&1&1&1&37&         131681894400&          81&         1&\text{  7}&\text{19}&\text{2A}&\text{4C}&\text{71}&\text{A5}&\text{CB}&\text{D6}\\
\hline\hline
\end{array}
\]
\renewcommand{\arraystretch}{1.5}
}
\end{table}

\renewcommand{\arraystretch}{0.3}
{\tiny
\begin{table}
\caption{Transposed incidence matrix $M_7$
containing all $M_n$,
$1\le n\le7$.
Symbol at position $(s',s)$
carries information about $m_{s,s'}$~\eqref{SNFinc}:
$\bullet$, $\star$, $\circ$ denotes respectively $1$,
$0$ explained by Lemma~\ref{SNFogr}, and
$0$ not explained by Lemma~\ref{SNFogr}.
}
\label{SNFincid}
\begin{tabular}{||r@{\,}r||@{\rg}c@{\rg}|@{\rg}c@{\rg}||@{\rg}c@{\rg}||@{\rg}c@{\rg}||@{\rg}c@{\rg}||@{\rg}c@{\rg}|@{\rg}c@{\rg}||@{\rg}c@{\rg}|@{\rg}c@{\rg}||@{\rg}c@{\rg}|@{\rg}c@{\rg}|@{\rg}c@{\rg}||}
\hline\hline
        &  & 0  & 0  & 0  & 00    & 000      & 00000          &0  & 000000000                  &000      & 111111111111111111                                    &11111          &1111111               \\
        &  & 0  & 0  & 0  & 00    & 000      & 00000          &0  & 111111111                  &111      & 111111111111111111                                    &11111          &1111111               \\
 ~~~$s$ &  & 0  & 0  & 0  & 00    & 000      & 11111          &1  & 111111111                  &111      & 111111111111111111                                    &11111          &1111111               \\
        &  & 0  & 0  & 0  & 00    & 111      & 11111          &1  & 111111111                  &111      & 111111111111111111                                    &11111          &1111112               \\
        &  & 0  & 0  & 0  & 11    & 111      & 11111          &1  & 111111111                  &112      & 111111111111111111                                    &11111          &1112222               \\
        &  & 0  & 0  & 1  & 11    & 111      & 11111          &2  & 111111111                  &222      & 111111111111111111                                    &22222          &3342222               \\
\cline{3-14}
$s'$~~~ &  & 0  & 0  & 0  & 00    & 000      & 00000          &0  & 000000000                  &000      & 000000000111111111                                    &00001          &0000000               \\
\hline\hline
\end{tabular}
\end{table}
}

\renewcommand{\arraystretch}{0.23}
{\tiny
\begin{table}
\addtocounter{table}{-1}
\caption{Continued.}
\centering
\begin{tabular}{||r@{\,}r||@{\rg}c@{\rg}|@{\rg}c@{\rg}||@{\rg}c@{\rg}||@{\rg}c@{\rg}||@{\rg}c@{\rg}||@{\rg}c@{\rg}|@{\rg}c@{\rg}||@{\rg}c@{\rg}|@{\rg}c@{\rg}||@{\rg}c@{\rg}|@{\rg}c@{\rg}|@{\rg}c@{\rg}||}
\hline\hline
        &  & 0  & 0  & 0  & 00    & 000      & 00000          &0  & 000000000                  &000      & 111111111111111111                                    &11111          &1111111               \\
        &  & 0  & 0  & 0  & 00    & 000      & 00000          &0  & 111111111                  &111      & 111111111111111111                                    &11111          &1111111               \\
 ~~~$s$ &  & 0  & 0  & 0  & 00    & 000      & 11111          &1  & 111111111                  &111      & 111111111111111111                                    &11111          &1111111               \\
        &  & 0  & 0  & 0  & 00    & 111      & 11111          &1  & 111111111                  &111      & 111111111111111111                                    &11111          &1111112               \\
        &  & 0  & 0  & 0  & 11    & 111      & 11111          &1  & 111111111                  &112      & 111111111111111111                                    &11111          &1112222               \\
        &  & 0  & 0  & 1  & 11    & 111      & 11111          &2  & 111111111                  &222      & 111111111111111111                                    &22222          &3342222               \\
\cline{3-14}
$s'$~~~ &  & 0  & 0  & 0  & 00    & 000      & 00000          &0  & 000000000                  &000      & 000000000111111111                                    &00001          &0000000               \\
\hline\hline
\end{tabular}
\end{table}
}\renewcommand{\arraystretch}{1.5}

\renewcommand{\arraystretch}{1.1}
\begin{table}
\caption{The part of $\mathcal{S}_9$ corresponding to nonsingular part of $\mathcal{A}_9$.}
\label{SNFdevet}
{\small
\begin{tabular}{||r|r|rrrrrrrrr|c||}
\hline\hline
        &    det &\multicolumn{10}{c||}{SNF} \\
\hline\hline
1-106   &    $m$ &   1&   1&   1&   1&   1&   1&   1&   1& $m$ &
$m\in\{1-98,100,101,102,\}$ \\
        &        &    &    &    &    &    &    &    &    &      &
$\cup\{104,105,108,110,120\}$ \\
\hline
107-135 &  $4m$  &   1&   1&   1&   1&   1&   1&   1&   2& $2m$ &
$m\in\{1-29\}$ \\
\hline
136-148 &    $9m$&   1&   1&   1&   1&   1&   1&   1&   3& $3m$ &
$m\in\{1-13\}$ \\
\hline
149-155 &   $16m$&   1&   1&   1&   1&   1&   1&   1&   4& $4m$ &
$m\in\{1-7\}$ \\
\hline
156-159 &      25&   1&   1&   1&   1&   1&   1&   1&   5& $5m$ &
$m\in\{1-3,5\}$ \\
\hline
 160    &      36&   1&   1&   1&   1&   1&   1&   1&   6&   6& \\
\hline\hline
\end{tabular}
}
\end{table}

\renewcommand{\arraystretch}{0.3}

\begin{table}
\caption{The part of the transposed incidence matrix $M_8$,
corresponding to the regular part of
$\mathcal{S}_{9}$.
Symbol at position $(s',s)$
carries information about $m_{s,s'}$~\eqref{SNFinc}:
$\bullet$, $\star$, $\circ$ denotes respectively $1$,
$0$ explained by Lemma~\ref{SNFogr}, and
$0$ not explained by Lemma~\ref{SNFogr}.
}
\label{SNFincidb}
{\tiny
\begin{tabular}{||@{\rx}r@{\,}r@{\,}r@{\rx}||@{\rx}c@{\rx}|@{\rx}c@{\rx}|@{\rx}c@{\rx}||@{\rx}c@{\rx}|@{\rx}c@{\rx}|@{\rx}c@{\rx}|@{\rx}c@{\rx}|@{\rx}c@{\rx}|@{\rx}c@{\rx}||}
\hline\hline
        &  &   & 000000000000000000                                    &0000        &00000000                &11111111111111111111111111111111111111111                                                                                  &111111111111                        &11111          &1111        &1111111              &1111         \\
        &  &   & 111111111111111111                                    &1111        &11111111                &11111111111111111111111111111111111111111                                                                                  &111111111111                        &11111          &1111        &1111111              &1111         \\
        &  &   & 111111111111111111                                    &1111        &11111111                &11111111111111111111111111111111111111111                                                                                  &111111111111                        &11111          &1111        &1111111              &1111         \\
 ~~~$s$ &  &   & 111111111111111111                                    &1111        &11111111                &11111111111111111111111111111111111111111                                                                                  &111111111111                        &11111          &1111        &1111111              &1111         \\
\hline\hline
\end{tabular}
}
\end{table}

\begin{table}
\addtocounter{table}{-1}
\caption{Continued.}
\centering
{\tiny
\begin{tabular}{||@{\rx}r@{\,}r@{\,}r@{\rx}||@{\rx}c@{\rx}|@{\rx}c@{\rx}|@{\rx}c@{\rx}||@{\rx}c@{\rx}|@{\rx}c@{\rx}|@{\rx}c@{\rx}|@{\rx}c@{\rx}|@{\rx}c@{\rx}|@{\rx}c@{\rx}||}
\hline\hline
        &  &   & 000000000000000000                                    &0000        &00000000                &11111111111111111111111111111111111111111                                                                                  &111111111111                        &11111          &1111        &1111111              &1111         \\
        &  &   & 111111111111111111                                    &1111        &11111111                &11111111111111111111111111111111111111111                                                                                  &111111111111                        &11111          &1111        &1111111              &1111         \\
        &  &   & 111111111111111111                                    &1111        &11111111                &11111111111111111111111111111111111111111                                                                                  &111111111111                        &11111          &1111        &1111111              &1111         \\
 ~~~$s$ &  &   & 111111111111111111                                    &1111        &11111111                &11111111111111111111111111111111111111111                                                                                  &111111111111                        &11111          &1111        &1111111              &1111         \\
\hline\hline
\end{tabular}
}
\end{table}

\begin{table}
\addtocounter{table}{-1}
\caption{Continued.}
\centering
{\tiny
\begin{tabular}{||@{\rx}r@{\,}r@{\,}r@{\rx}||@{\rx}c@{\rx}|@{\rx}c@{\rx}|@{\rx}c@{\rx}||@{\rx}c@{\rx}|@{\rx}c@{\rx}|@{\rx}c@{\rx}|@{\rx}c@{\rx}|@{\rx}c@{\rx}|@{\rx}c@{\rx}||}
\hline\hline
        &  &   & 000000000000000000                                    &0000        &00000000                &11111111111111111111111111111111111111111                                                                                  &111111111111                        &11111          &1111        &1111111              &1111         \\
        &  &   & 111111111111111111                                    &1111        &11111111                &11111111111111111111111111111111111111111                                                                                  &111111111111                        &11111          &1111        &1111111              &1111         \\
        &  &   & 111111111111111111                                    &1111        &11111111                &11111111111111111111111111111111111111111                                                                                  &111111111111                        &11111          &1111        &1111111              &1111         \\
 ~~~$s$ &  &   & 111111111111111111                                    &1111        &11111111                &11111111111111111111111111111111111111111                                                                                  &111111111111                        &11111          &1111        &1111111              &1111         \\
\hline\hline
\end{tabular}
}
\end{table}

\renewcommand{\arraystretch}{1.5}

\begin{table}
\caption{Lower bounds for $a_n$
and matrices with high extension spectra, $10\le n\le 19$.}
\label{prve0}
{\small
\begin{tabular}{||r|r||r|r|r|l||}
\hline\hline
$n$ & $a_n\ge$ & $|\mathcal{L}_{n-1}|$ & $\det A_{n-1}$ & $a_{A_{n-1}}$  & $A_{n-1}$   \rule[-2mm]{0cm}{6mm}   \\
\hline\hline
 $ 10$&$    259$&$       2 $&$      110$&$     257$&$  [ 7,39,5A,9C,E1,149,174,193,1AA]$\\
 $ 11$&$    739$&$       6 $&$      291$&$     679$&$  [ F,71,B6,13A,1C3,1DC,256,299,2EC,325] $\\
 $ 12$&$   2107$&$      19 $&$      779$&$    1894$&$  [ F,73,195,1EA,2A6,35C,4D6,53E,565, $\\
      &         &           &           &          &    $ 6B9,703]$               \\
 $ 13$&$   6157$&$      18 $&$     2201$&$    5618$&$  [ 1F,E3,17C,3A5,649,6D6,732,A6E,AB8,$\\
      &         &           &           &          &    $ B53,C35,D8E] $                       \\
 $ 14$&$  19073$&$      40 $&$     6731$&$   16821$&$  [ 3F,1C7,2D9,76A,C4D,CF2,F94,1575,$\\
\hline\hline
\end{tabular}
}
\end{table}

\end{document}